\documentclass[reqno]{amsart}
\usepackage{graphicx} 
\usepackage{amsmath, amssymb ,amsthm, amsfonts, amsgen, color, mathtools}
\usepackage{graphicx, hyperref, esint}
\usepackage{bbm}
\usepackage{tikz}
\usepackage[colorinlistoftodos]{todonotes}
\usepackage{hyperref}
\usepackage{mathrsfs}
\usepackage[OT1]{fontenc}
\usepackage[bottom]{footmisc}
\usepackage{mathabx}
\numberwithin{equation}{section}
\definecolor{vg}{rgb}{0.0, 0.26, 0.15}

\newcommand{\R}{{\mathbb{R}}}

\newcommand{\A}{{\mathcal{A}}}
\newcommand{\cL}{{\mathcal{L}}}
\newcommand{\beq}{\begin{equation}}
\newcommand{\eeq}{\end{equation}}

\newcommand{\Mcal}{{\mathcal{M}}}

\newcommand{\eps}{{\epsilon}}

\newcommand{\wsto}{\stackrel{*}{\rightharpoonup}}

\renewcommand{\eps}{\varepsilon}

\def\e{{\varepsilon}}

\def\A{\mathcal A}
\def\B{\mathcal B}

\def\Lin{\mbox{Lin}}

\def\symB{\mathbb{B}}

\newcommand{\average}{{\mathchoice {\kern1ex\vcenter{\hrule
				height.4pt width 8pt depth0pt}
			\kern-11pt} {\kern1ex\vcenter{\hrule height.4pt width 4.3pt
				depth0pt} \kern-7pt} {} {} }}
\newcommand{\ave}{\average\int}

\def\XXint#1#2#3{{\setbox 0=\hbox{$#1{#2#3}{\int}$}
\vcenter{\hbox{$#2#3$}}\kern-.5\wd0}}

\definecolor{vg}{rgb}{0.0, 0.26, 0.15}
\definecolor{coolblack}{rgb}{0.0, 0.18, 0.39}

\newcommand{\res}{\mathop{\hbox{\vrule height 7pt width .5pt depth 0pt
\vrule height .5pt width 6pt depth 0pt}}\nolimits}

\newtheorem{Theorem}{Theorem}[section]
\newtheorem{lemma}[Theorem]{Lemma}
\newtheorem{Proposition}[Theorem]{Proposition}
\newtheorem{Corollary}[Theorem]{Corollary}

\newtheorem{Remark}[Theorem]{Remark}
\newtheorem{Definition}[Theorem]{Definition}
\newtheorem{Examples}[Theorem]{Examples}
\def\f1min{f_1^{\rm min}}

\usepackage{color}
\usepackage{natbib}
\setcitestyle{numbers,square,semicolon}

\begin{document}
\title{Towards a unified treatment of structured deformations theories}

\author[L. D'Elia]{Lorenza D'Elia}
\address[L. D'Elia]{Institute of Analysis and Scientific Computing, TU Wien, Wiedner Hauptstraße 8-10, 1040 Vienna, Austria \\ \textit{ORCID}: \href{https://orcid.org/0000-0002-2809-5553}{0000-0002-2809-5553}}
\email{lorenza.delia@tuwien.ac.at}
  \par

\author[E.~Zappale]{Elvira Zappale}
\address[E.~Zappale]{Dipartimento di Scienze di Base ed applicate per l'Ingegneria, Sapienza - Universit\`{a} di Roma, Via A. Scarpa, 16, 00161, Roma, Italy, \& Research Centre in Mathematics and Applications,
University of \'{E}vora, Portugal \\ \textit{ORCID}: \href{https://orcid.org/0000-0001-7419-300X}{0000-0001-7419-300X}}
\email{elvira.zappale@uniroma1.it}

	

\begin{abstract}
Relying on the theory of $\mathcal B$ gradients, with $\mathcal B$ a $\mathbb C$-elliptic homogeneous linear partial differential operator of order $k$, with constant coefficients, we propose a method which allows to treat in a unified way the theories of structured deformations. 
Introducing the space of structured $\mathcal B$-deformations, $SD_{\mathcal B}$, we establish an integral-representation result for the associated energies, involving both bulk and surface contributions. Our approach relies on a direct application of the global method for relaxation, suitably adapted to the space of functions of $\mathcal B$-bounded variation. 
\color{black}
\end{abstract}

\maketitle

\section{Introduction}
The theory of \emph{structured deformations} $SD$, introduced in a seminal work  by Del Piero and Owen \cite{DPO}, captures the effects at the macroscopic level of both smooth and non\--\-smooth deformations at a sub-macroscopic level. More precisely, the deformation of a body is characterized by a pair of fields $(g,G)$ and a set $\mathcal K$  (a surface-like subset of the reference body $\Omega \subset \mathbb R^N$) called a \emph{structured deformation}: the macroscopic deformation field is $g$ and $G$ is the contribution at the macroscopic scale of smooth sub-macroscopic changes, while the difference $\nabla g-G$ captures the contribution at the macroscopic scale of non-smooth sub-macroscopic changes, such as slips and separations (referred to as \emph{disarrangements} \cite{DO2002}),  $\mathcal K$ plays the role of the ``crack site'', i.e. $g$ is regular in $\Omega \setminus \mathcal K$.

The consistency of Del Piero and Owen's $SD$ theory is ensured by the so-called {\it approximation theorem} stating that each structured deformation  $(g, G, \mathcal K)$ can always be approximated in an appropriate sense by a sequence $(g_n, \mathcal K_n)_n$ of classical deformations. 

Choksi and Fonseca \cite{CF1997} rephrased the theory in variational terms, namely setting the 'classical' part $g$ and the set  $\mathcal K$ as a single element $g$ of $SBV(\Omega;\mathbb R^d)$ (see \cite{AFP}), via the identification of $\mathcal K$ with the discontinuity set $J_g$. Therefore, in \cite{CF1997} a structured deformation is defined as a pair $(g,G)\in SBV(\Omega;\R^{d})\times L^1(\Omega;\R^{d\times N})$. 
Moreover, an {\it approximation theorem} with respect to a suitable notion of convergence has been proved. Thus, departing from an energy featuring a bulk contribution measuring  the deformation (gradient) and  a surface  contribution accounting for the energy needed to fracture the body, i.e. 
a functional of the type 
\begin{equation*}
	F(u)\coloneqq \int_\Omega \mathscr{W}\big(x,\mathcal \nabla u(x)\big)\,d x+\int_{\Omega\cap J_u} \psi\big(x,[u](x), \nu_u(x)\big)\,d \mathcal H^{N-1}(x)
\end{equation*}
 for  fields   $u\in SBV(\Omega;\mathbb{R}^d)$,  with
$\mathscr{W}\colon\Omega\times  \mathbb R^{d\times N}  \to[0,+\infty)$, $\psi\colon\Omega\times \mathbb R^{d} \times\mathbb S^{N-1}\to[0,+\infty)$, $\nu_u$ and $[u]$ being the unit normal to the crack site and the size of the crack respectively,  the energy associated to a structured deformation is defined as the ``most economical way'' to obtain it approaching the given structured deformation by classical ones. More precisely, the energetic functionals are defined as 
\begin{equation}\label{102intro}
	I_{SD}(g,G)\coloneqq \inf\Big\{\liminf_{n\to\infty} F(u_n)\colon \,  u_n \in SBV(\Omega\;\mathbb R^d), \ u_n\leadsto (g,G)\Big\},
\end{equation}
 where  $\leadsto$ denotes the convergence justified by the {\it approximation theorem} proved in \cite{DPO} and \cite{CF1997}.

In more recent years  the theory of {\it first-order} structured deformations has been broadened also in view of its applications \cite{DO2003, PDOF}, see \cite{MMO} and the references therein for a  more comprehensive survey about this theory  which include, among others, optimal design, dimension reduction, and homogenization. 
Subsequently, Owen and Paroni \cite{OP} extended the theory to {\it second\--\-order} structured deformations, and \cite{O17} presented models accounting for strain-gradient theories which allow other geometrical changes at a sub-macroscopic level, such as curvature, or gradient-disarrangements. 

While the variational formulation associated to \cite{OP} was carried out by \cite{FLP, FHP}(see also \cite{BMMO2017} for another model which features  jumps for both the approximating fields as well as for their gradients), the strain-gradient setting discussed in \cite{O17} has not be addressed so far in the literature. This paper, in fact, proposes a first variational modeling in this framework, see Examples \ref{examples}.  

In \cite{FMZ,MSZ} and \cite{DF} the theory of structured deformations has been developed in the context of linearized elasticity where, under the assumption of infinitesimally small deformations, the material response depends only the symmetric part of the displacement gradient.
In this setting the space where a structured deformation $(g, G)$ lies is $SBD\times L^1$, denoted in \cite{FMZ} as $StSBD\times L^1$, which in our current notations is $SD_{\mathscr E}$, where the generic differential operator $\B $ is specialized in the distributional symmetrized gradient $\mathscr{E}$.
The motivation for such an analysis lies in the application to solid mechanics. Indeed, if deformation gradients are very close to the set of orientation preserving rigid motions, the theory of linear elasticity is expected to deliver a suitable description. 

It is clear that to deal with different materials' behaviours, it would be interesting to consider also more general deformation gradient theories and accounting for possible submacroscopic disarrangements. To this end, we address here the case of deformations $g \in BV^\B(\Omega)$, with $\mathcal B$  not necessarily $\mathscr{E}, D$ or $D^2$ and $G$ a summable field defined in $\Omega$ with values in ${Im B}$, see Subsection \ref{subsec:pseudoop} and Section \ref{sec:functionspaces} for precise definitions. In order to fix ideas, let  $\Omega$ be an open, bounded subset of $\mathbb{R}^N$ and let $\B$ be a $\mathbb C$-elliptic homogeneous linear partial differential operator of order $k\geq 1$, between finite-dimensional (inner product) Euclidean spaces $V $ and $W$. More precisely, we shall consider operators acting on smooth maps $u \in C^\infty(\Omega; V)$ as
\begin{equation}\label{Bdef}
\mathcal Bu = \sum_{|\beta| = k} B_\beta \, \partial^\beta u \in C^\infty(\Omega; W),
\end{equation}
otherwise the action on summable fields $u$ has to be intended distributionally. Here, the coefficients $B_\beta \in \mathrm{Lin}(V; W)$ are assumed to be constant, $\beta \in \mathbb{N}_0^n$ is a multi-index with modulus $|\beta| = \beta_1 + \cdots + \beta_n$, and $\partial^\beta$ denotes the distributional derivative $\partial_1^{\beta_1} \cdots \partial_n^{\beta_n}$. Such operator $\B$ can be expressed through the \emph{jet-form}, which is defined as 
\begin{align}\label{jetformBu}
B^k(D^k u) = \mathcal Bu,
\end{align}
where $B^k : V \otimes E_k(\mathbb{R}^N) \to W$ is the (unique) linear map satisfying
\[
B^k[v \otimes \xi^{\otimes k}] = {\mathbb B}^k(\xi)[v]
\quad \text{for all } v \in V,
\notag
\]
with ${\mathbb B}$ is the (principal) symbol associated to $\B$ (see Subections \ref{sec:Notation} and \ref{subsec:pseudoop} for the definitions and the properties of $\mathbb B$ and  $E_k(\R^N)$). In the following, we will denote by ${\rm Im}B^k$  the vector subspace of $W$ which is the  image of $V\otimes E_k(\mathbb R^N)$ through $B^k$, namely
\begin{align}\label{ImBk}
{\rm Im}B^k:=\{ \eta \in W: \exists \zeta \in V\otimes E_k(\mathbb R^N) \hbox{ such that } \eta = B^k \zeta\}.
\end{align}

To encompass the above mentioned existing theories of structured deformations and, at the same time, to have the possibility of dealing with different material responses, we introduce the class of $k$th-order structured deformations
$SD_{\B}(\Omega)$, as follows.

\begin{Definition}\label{def:SD}
The space of $k$th-order structured deformations
$SD_{\B}(\Omega)$, relative to the $\mathbb C$-elliptic, differential homogeneous operator $\B$, of $k$-th order, consists of pairs $(u,U)$ with
$u\in SBV^{\B}(\Omega)$ and
$U\in L^1(\Omega; {\rm Im}B^k)$, i.e., 
\[
SD_{\B}(\Omega)
:=
SBV^{\B}(\Omega)
\times
L^1(\Omega; {\rm Im}B^k).
\]
\end{Definition}    
We direct the reader to Section \ref{sec:functionspaces} for the introduction of the space $SBV^\B(\Omega)$ of special fields with bounded $\B$-variation.  It suffices, now, thinking of the space $SBV^\B(\Omega)$ as the space of Sobolev  fields $u \in  W^{k-1,1}(\Omega;V)$ whose distributional $k$-th order $\B$-gradient $\mathcal B u$ is a bounded Radon measure decomposable in a bulk part, absolutely continuous with respect to the Lebesgue measure and surface type measure, concentrated where the $\nabla^{k-1} u$ has different traces, see Definition \ref{def:approxjumppoint}.    

Thus, in perfect analogy with the { \it approximation theorem} in \cite{CF1997}, we  prove in Section $4$ the following result, when $\B$ is a $\mathbb C$-elliptic operator as in \eqref{Bdef} and  of order $k$ at most $2$. The restriction on the order of the operator is due to the possibility of interpreting any field in $\mathbb R^d \otimes E_k(\mathbb R^N)$ (here $V \cong \mathbb R^d$) as the absolutely continuous part of $k$-th order gradient, see Alberti's Theorem  \ref{Al} below for the case $k=1$.

\begin{Theorem}\label{appthintro}
 Let $\mathcal{B}$ be a $\mathbb C$-elliptic linear homogeneous partial differential operator of order $k\leq 2$, whose symbol in jet form is $B^k$.  Let ${\rm Im}B^k\subseteq W$ be as in \eqref{ImBk}. Let  $(g,G)\in SBV^{\B}(\Omega)\times L^1(\Omega;{\rm Im} B^k)$.
 Then,  there exists a sequence $(g_n)_n \subset SBV^{\B}(\Omega)$ such that
$g_n \wsto g$ in $BV^{\mathcal B}(\Omega)$ and $\mathcal B^a g_n \wsto G$ in $\mathcal M(\Omega;{\rm Im}B^k)$.
\end{Theorem}

Therefore, to any deformation in the space $SBV^\B(\Omega)$, where different gradient theories are encompassed (see \cite{B} and Examples \ref{examples}), it is possible to associate an energy, which is the sum of a bulk contribution measuring the deformation (that in our case $\B$-gradient) and  a surface  contribution accounting for the energy needed to fracture the body, i.e. 
a functional of the type 
\begin{equation*}
	F_{\B}(u)\coloneqq \int_\Omega \mathscr{W}\big(x,\B u(x)\big)\,d x+\int_{\Omega\cap J_u} \psi\big(x,[u](x), \nu_u(x)\big)\,d \mathcal H^{N-1}(x)
\end{equation*}
 for  fields   $u\in SBV^\B(\Omega)$,  with
$\mathscr{W}\colon\Omega\times W  \to[0,+\infty)$ and $\psi\colon\Omega\times \mathbb R^{d} \times\mathbb S^{N-1}\to[0,+\infty)$. Then, the energy associated to a $\B$-structured deformation is defined as in \cite{CF1997}, via a relaxation: 
\begin{equation}\label{102introB}
	I_{SD_\B}(g,G)\coloneqq \inf\Big\{\liminf_{n\to\infty} F_\B(u_n)\colon \,  u_n \in SBV^\B(\Omega), \ u_n\leadsto (g,G)\Big\},
\end{equation}
 where  $\leadsto\equiv (\wsto,\wsto)$ denotes the convergence in Theorem \ref{appthintro}.
 To achieve an integral representation for $I_{SD_\B}$, we prove in Section \ref{gm} a general integral representation result, i.e. a {\it Global Method} for relaxation in $SD_\B$. Our method is tailored  for structured deformations in the  full $\mathcal B$-gradients setting with the additional geometrical field $G$ belonging to $L^1$ and having the symmetries inherited by the operator $\B$ via the symbol $B^k$, see  Theorem \ref{thm:GMSDB}.  In Section \ref{gm}, we actually obtain a generalization of the global method for relaxation, introduced by Bouchitt\'e, Fonseca, and Mascarenhas \cite{BFM} in $BV$-setting, and later addressed in the many settings (see, for a not exhaustive list \cite{ARS-24}, \cite{CFVG}, \cite{SSS}, \cite{FMZ} and the references therein). In fact, we provide in Theorem \ref{thm:globalmethod} a general method for the identification of the integral representation at least on $SBV^\B$ of a class of
functionals defined in  $BV^\B$.   
This is indeed sufficient for further applications to relaxation in different contexts, such as the integral representation of specific functionals in $BV^\B\times \mathcal M$ in \cite[Section 4]{DEZ}, see Remark \ref{rem:caseBVxL1}. 
  
Let $\mathcal A(\Omega)$ be the class of open subsets  of $\Omega$, we consider a functional
\begin{equation}\label{def:introF}
\mathcal{F}:
SD_{\B}(\Omega)\times\mathcal{A}(\Omega)
\longrightarrow [0,+\infty)
\end{equation}
satisfying the following hypotheses:

\begin{itemize}
\item[($\tilde{H}$1)]
$\mathcal{F}(u,U;\cdot)$ is the restriction to
$\mathcal{A}(\Omega)$ of a Radon measure for every
$(u,U)\in SD_{\B}(\Omega)$.

\item[($\tilde{H}$2)]
$\mathcal{F}(\cdot,\cdot;A)$ is
$SD_{\B}$-lower semicontinuous, in the sense that if
$(u,U)\in SD_{\B}(\Omega)$,
$(u_n,U_n)\subset SD_{\B}(\Omega)$,
$u_n\to u$ in $L^{1}(\Omega;\mathbb{R}^d)$
and $U_n\stackrel{*}{\rightharpoonup}U$
in $\mathcal{M}(\Omega)$, then
\[
\mathcal{F}(u,U;A)
\leq
\liminf_{n\to+\infty}
\mathcal{F}(u_n,U_n;A).
\]

\item[($\tilde{H}$3)]
$\mathcal{F}$ is local, i.e., for all
$A\in\mathcal{A}(\Omega)$, if $u=v$ and $U=V$
$\mathcal{L}^N$-a.e. $x\in A$, then
\[
\mathcal{F}(u,U;A)=\mathcal{F}(v,V;A).
\]

\item[($\tilde{H}$4)]
There exists a constant $C>0$ such that
\[
\frac{1}{C}
\left(
\|U\|_{L^1(A)}+|\B u|(A)
\right)
\leq
\mathcal{F}(u,U;A)
\leq
C\left(
\mathcal{L}^N(A)
+\|U\|_{L^1(A)}
+|\B u|(A)
\right)
\]
for every $(u,U)\in SD_{\B}(\Omega)$,
$A\in\mathcal{A}(\Omega)$.
\end{itemize}

In the spirit of the global method for relaxation,
Bouchitt\'e et al. (1998, 2002), given
$(u,U;A)\in SD_{\B}(\Omega)\times\mathcal{A}(\Omega)$
we define
\begin{equation}\label{compSDB}
\mathfrak{A}(u,U;A)
:=
\left\{
(v,V)\in SD_{\B}(\Omega):
\operatorname{spt}(u-v)\Subset A,\,
\int_A(U-V)\,dx=0
\right\},
\end{equation}
and
\begin{equation}
m(u,U;A)
:=
\inf\left\{
\mathcal{F}(v,V;A):
(v,V)\in\mathfrak{A}(u,U;A)
\right\}.
\label{m4.3}
\end{equation}

\begin{Theorem}\label{thm:GMSDB}
Assume Assumptions ($\tilde H$1), ($\tilde H$2), ($\tilde H$3) and ($\tilde H$4). Then, 
for every $(u,U)\in SD_{\B}(\Omega)$ and
$A\in\mathcal{A}(\Omega)$, we have
\begin{equation}\notag
\begin{aligned}
\mathcal{F}(u,U;A)
&=
\int_A f(x,u,\nabla u,\dots, \nabla^{k} u,U)\,dx\\
&\quad+
\int_{S(\nabla^{k-1} u)\cap A}
h(x,u,\nabla u, \dots, \nabla^{k-2} u, (\nabla^{k-1} u)^+,(\nabla^{k-1} u)^-,\nu_{\nabla^{k-1} u})
\,d\mathcal{H}^{N-1},
\end{aligned}
\end{equation}
where 
  \begin{equation}
      \label{fdef-U}
      f(x_0, G_0, G_1, \dots, G_k, H) = \limsup_{\e\to 0} \frac{m(P^k_{x_0}(x-x_0), H; Q(x_0, \e))}{\e^N},
  \end{equation}
where for every $j \in \{1,\dots, k\}$, $P^j_{x_0}(x-x_0)$ is the $j$-th order polynomial of the form
             \begin{equation}
                 \label{Pkx0def}
                 P^j_{x_0}(x-x_0)= G_0+ \langle G_1,x-x_0\rangle + \tfrac{1}{2} \langle G_{2}, (x-x_0)^{\otimes^2} \rangle +\dots
                 + \tfrac{1}{j!} \langle G_j, (x-x_0)^{\otimes ^j}\rangle,   
             \end{equation}
             with $G_j \in \mathbb V \otimes E_j(\mathbb R^N)$, for $j=0,\dots, k$, and 
\begin{align}\label{hdef-U}
h(x_0, G_0, G_1,\dots, G_{k-2}, c, d, \nu)=\limsup_{\varepsilon \to 0^+} \frac{m(P_{x_0}^{k-2}(x-x_0)+v_{J}(x), O; Q_\nu(x_0, \e))}{\varepsilon^{N-1}}
\end{align}
where $O\in\mathbb{R}^{d\times N^k}$ is the matrix with all entries equal to zero and $v_J$ is defined as 
\begin{align}\label{v_J}
 v_J(x) \coloneqq v_{c,d,\nu}(x-x_0):=\frac{1}{(k-1)!}\left\{
\begin{array}{ll}
\langle c, (x-x_0)^{\otimes^{k-1}}\rangle &\hbox{ if } (x-x_0)\cdot \nu > 0,\\
\langle d, (x-x_0)^{\otimes^{k-1}}\rangle &\hbox{ if }(x-x_0)\cdot \nu \leq 0. 
\end{array}\right.
\end{align}
\end{Theorem}

The proof of Theorem \ref{thm:GMSDB} hinges on a Poincar\'{e} inequality (see Corollary \ref{cor:poincare}) in $BV^{\B}$. The latter follows from an application of a general result of Ziemer \cite[Lemma 4.1.3]{Z} (see also Lemma \ref{lemma:Ziemerlemma4.1.3}) together with the compactness result proved in \cite[Corollary 2.16]{DEZ}. 

It is worth to emphasize that in particular cases, such as $\mathcal B=D$ or $\mathscr E$, Theorem \ref{thm:GMSDB} can be applied to several relaxation problems in the contexts of $SD_\B$, obtaining a more detailed description of the energy densities representing $I_{SD_\B}$. In particular one can recover, among the others, the results in \cite{CF1997, FHP, AMMZ, BMZ, FMZ}. Moreover when there is enough information on the operator $\B$ (so far partial results are available in the general case, see \cite{AR, C26}), and Federer-Vol'pert type-theorem (see, e.g. in the $BV$ case \cite[Theorem 3.78]{AFP}) are available, via Remark \ref{rem:caseBVxL1}, the theory of $\mathcal B$-structured deformations can comprise the entire $BV^\B$ allowing to recover the results in \cite{KKMZ, FMZ}. 
\medskip

The paper is organized as follows: after introducing the problem and setting the notation, in Section \ref{sec:preliminaries} we provide some preliminary results that will be useful in the subsequent  sections.  Section \ref{gm} is dedicated to abstract integral representation results, i.e. to a global method in the space of fields of generalized $\mathcal B$ bounded variation. Section \ref{sec:applications} concerns not only the variational theory of $\mathcal B$-structured deformations, including approximation theorems, relaxation and integral representation results, but some application of the global method to the space $BV^\B(\Omega)$.
Finally, in the Appendix \ref{sec:appendix} we present some results about measure theory. We also provide some results dealing with approximate differentiability in the spaces $BV^\B$.

\medskip

\section{Preliminaries}\label{sec:preliminaries}
\subsection{Notation} 
\label{sec:Notation}
Throughout the present manuscript, we will use the following notations
\begin{itemize}
	\item $\mathbb N$ denotes the set of natural numbers without the zero element;
	\item $\mathbb S^{N-1}$ denotes the unit sphere in $\mathbb R^N$;
    \item $V, W$ are vector spaces of finite dimensions on $\mathbb K \coloneqq \{\mathbb R, \mathbb C\}$;
	\item $Q\coloneqq (-\tfrac12,\tfrac12)^N$ denotes the open unit cube of $\mathbb R^{N}$ centred at the origin; for any $\nu\in\mathbb S^{N-1}$, $Q_\nu$ denotes any open unit cube in $\mathbb R^{N}$ with two faces orthogonal to $\nu$; 
    \item $B_r(x)$ is the open unit ball of radius $r>0$ and centered at $x\in\R^N$;
	\item for any $x\in\mathbb R^{N}$ and $\delta>0$, $Q(x,\delta)\coloneqq x+\delta Q$ denotes the open cube in $\mathbb R^{N}$ centred at $x$ with side length $\delta$; likewise $Q_\nu(x,\delta)\coloneqq x+\delta Q_\nu$;
	\item ${\mathcal A}(\Omega)$ is the family of all open subsets of $\Omega $, whereas ${\mathcal A}_\infty(\Omega)$ is the family of all open subsets of $\Omega $ with Lipschitz boundary; 
	\item $\mathcal L^{N}$ and $\mathcal H^{N-1}$ denote the  $N$-dimensional Lebesgue measure and the $\left(  N-1\right)$-dimensional Hausdorff measure in $\mathbb R^N$, respectively; the symbol $\partial x$ will also be used to denote integration with respect to $\mathcal L^{N}$; 
	\item $\mathcal M(\Omega;\mathbb R^{d\times N})$ is the set of finite matrix-valued Radon measures on $\Omega$; $\mathcal M ^+(\Omega)$ is the set of non-negative finite Radon measures on $\Omega$;
	given $\mu\in\mathcal M(\Omega;\mathbb R^{d\times N})$,  
	the measure $|\mu|\in\mathcal M^+(\Omega)$ 
	denotes the total variation of $\mu$;
	\item A sequence of measures $(\mu_n)_n$ converges \emph{area-strictly} to $\mu$ if
\begin{equation}\label{def:areaconv}
\mu_n \rightharpoonup \mu \quad \text{and} \quad
\int_\Omega \sqrt{1 + \left|\frac{d\mu_n}{d\mathcal{L}^N}\right|^2}\,dx + |\mu_n^s|(\Omega)
\to
\int_\Omega \sqrt{1 + \left|\frac{d\mu}{d\mathcal{L}^N}\right|^2}\,dx + |\mu^s|(\Omega).
\end{equation}
    \item Let $r \geq 0$ be an integer. For a vector $e \in \mathbb{R}^N$, we write $e^{\otimes r}$ to denote the $r$-fold tensor that results by taking the tensorial product of $e$ with itself $r$ times (with the convention $e^{\otimes 0} = 1$). Notice that
\[
E_r(\mathbb{R}^N) = \mathrm{span} \{ e^{\otimes r} : e \in \mathbb{R}^N \}.
\]

      \item We define a linear contraction $\langle\cdot,\cdot\rangle_r : V \otimes E_r(\mathbb{R}^N) \times E_r(\mathbb{R}^N) \to V$  (with the Frobenius inner product $(\cdot : \cdot)$ on $\otimes^r \mathbb{R}^N$) as
\[
\langle v \otimes M, e\rangle_r = v \, (M : e^{\otimes r})
\quad \text{for all } v \in V, \, M \in E_r(\mathbb{R}^N) \text{ and } e \in \mathbb{R}^N.
\]

Notice that under our convention for $r = 0$, we have $(v, \cdot)_0 = v$.

\end{itemize}

\subsection{Pseudodifferential operators}
\label{subsec:pseudoop}

We consider $k$-th order homogeneous linear partial differential operators $\B$ of the form \eqref{Bdef}. We now recall some important structural properties of $\B$.  \\
Let $\mathbb{K}=\{\R, \mathbb C\}$. An operator $\B$ as in \eqref{Bdef} is said to be of {\it $\mathbb K$-constant rank}  (see \cite{FM99, M3}) if  there exists an integer $r>0$ such that
\begin{align}\notag
{\rm rank}_{\mathbb K} \mathbb B(\xi ) = r\,\,\, \hbox{ for all } \xi \in  \mathbb K^N \setminus\{0\},
\end{align}
where the tensor-valued $k$-homogeneous polynomial
\begin{align}\label{mathbbAdef}
\mathbb B(\xi ) := \sum_{|\beta|=k}
\xi^\beta B_\beta   \in   \Lin(V; W),\;\;\; \xi  \in \mathbb K^N,
\end{align}
is the principal symbol associated to the operator $\mathcal B$. Here, $\xi^\beta  := \xi_1^{\beta_1}\dots \xi_N^{\beta_N}$.\\
An operator $\B$  of the form \eqref{Bdef} is {\it $\mathbb K$-elliptic} when the (principal) symbol map
\[
\mathbb B^k(\xi) := \sum_{|\beta| = k} \xi^\beta B_\beta \in \mathrm{Lin}(\mathbb{K} \otimes V; \mathbb{K} \otimes W), 
\quad \xi \in \mathbb{K}^n,
\]
is an injective linear map. Namely, there exists a positive constant $c$ such that
\[
|\mathbb B^k(\xi)v| \geq c |\xi|^k |v| 
\quad \text{for all } \xi \in \mathbb{K}^n \text{ and all } v \in \mathbb{K} \otimes V.
\]
Note that since $\mathcal B$ has  constant rank if ${\rm dim}({\rm Im}(\mathbb B(\xi)))$ is a constant $r$ independent of $\xi \in  \mathbb K^N \setminus \{0\}$, it follows that $\mathbb C$-elliptic operators $\B$ satisfies the constant rank properties by the Rank-Nullity Theorem.  

Finally, the operator $\B$ can be express in the jet-form $B^k$ given by \eqref{jetformBu} and recalled below
\begin{align}\notag
B^k(D^k u) = \mathcal Bu,
\end{align}
where $B^k : V \otimes E_k(\mathbb{R}^N) \to W$ is the (unique) linear map satisfying
\[
B^k[v \otimes \xi^{\otimes k}] = {\mathbb B}^k(\xi)[v]
\quad \text{for all } v \in V,
\notag
\]
with ${\mathbb B}$ being given in \eqref{mathbbAdef}. Note that the existence of $B^k$ follow from the fact that the principal symbol ${\mathbb B}$ is a $k$-homogeneous polynomial together with the universal property of the  tensor product. When $k=1$, we simply write $B$ instead of $B^1$.

\subsection{Function Spaces} 
\label{sec:functionspaces}
Let $\B$ be a $k$-th order homogeneous linear partial differential operator. Let $\Omega$ be an open and bounded subset of $\R^N$ with Lipschitz boundary.  For $1\leq p\leq \infty$, the Sobolev space $W^{\B, p}(\Omega)$ associated to the operator $\B$ is defined by
     \begin{equation}
         \notag
         W^{\B, p}(\Omega) := \left\{ u\in L^p(\Omega; V) : \B u\in L^p(\Omega; W)   \right\}.
     \end{equation}
It is a Banach space endowed with the norm
        \begin{equation}
            \notag
            \|u\|_{W^{\B, p}(\Omega)}:= \|u\|_{L^p(\Omega; V)} + \|\B u\|_{L^p(\Omega; W)}.
        \end{equation}
In the case $\B = \nabla^k$,  we simply write $W^{k, p}(\Omega; \R^m)$, where the latter denotes the classical Sobolev space, i.e., the space of functions whose distributional derivatives up to the order $k$ are $p$-summable. When $k=0$, we mean the $L^p(\Omega;V)$ space. 
For further details on these spaces see \cite{BDG, GR19}.
We also introduce the space $BV^\B(\Omega)$ of functions of bounded $\B$-variations defined by
      \begin{equation}\notag
          BV^\B (\Omega):=\{u\in L^1(\Omega; V) : \B u\in \Mcal(\Omega; W) \}.
      \end{equation}   
It is a Banach space endowed with the following norm
      \begin{equation}
          \label{def:normBVB}
          \|u\|_{BV^\B}:= \|u\|_{L^1(\Omega; V)} + |\B u|(\Omega),
      \end{equation}
where $|\cdot|:\mathcal M(\Omega;W)\to [0,+\infty)$ stands for the total variation of a measure.
If $\B$ is a $k$-th order $\mathbb C$-elliptic differential operator, the space $BV^\B$ can be seen as a generalization of $\nabla^k B$ spaces introduced in \cite[eq. (2.45) and (2.46)]{DT} (see also the spaces $BV^k$ introduced in \cite{ADC}, and recall that, when $k=2$, $\nabla^2 B= BV^2 = BH$ where $BH$ is the space of functions with bounded hessian, see \cite{D}). 
Thanks to \cite[Theorem 3.7]{DG24} (see also \cite[Theorem 2.2]{PPRV}), $\mathbb C$-ellipticity of $\B$ guarantees the following equality
\begin{align}\label{BVBequiv}BV^\B(\Omega)= W^{k-1, 1}(\Omega; V) \cap BV^\B(\Omega).
\end{align}

To introduce the jump set of functions $u\in BV^{\mathcal B}(\Omega)$, we only consider $\mathbb C$-elliptic operator $\mathcal B$ of first order. This is not restricted, since the order-reduction procedure developed in \cite{ARS} allows higher-order differential operators to be represented through suitable first-order operators. Consequently, jump-set properties established for first-order operators can be extended to operators of arbitrary order. We first recall the definition of \emph{approximate jump}.

\begin{Definition} \label{def:approxjumppoint}
Let $u\in L^1_{\mathrm{loc}}(\Omega; V)$.
A point $x$ is an \emph{approximate jump point}
of $u$, and we write $x\in J_u$, if there exist
distinct vectors $a,b\in V$ and a direction
$\nu\in S^{n-1}$ satisfying
\begin{equation}\label{eq:approxjumppoint}
\begin{cases}
\displaystyle
\lim_{r\downarrow 0}
\fint_{B_r^+(x,\nu)}
|u(y)-a|\,dy=0,
\\[1.5em]
\displaystyle
\lim_{r\downarrow 0}
\fint_{B_r^-(x,\nu)}
|u(y)-b|\,dy=0,
\end{cases}
\end{equation}
where the $\nu$-oriented half-balls $B_r^+(x,\nu)$ and $B_r^-(x,\nu)$ centred at $x$ are defined as 
\begin{equation} \label{def:halfball}
  B_r^+(x,\nu)
:=
\left\{
y\in B_r(x):
\langle\nu,y\rangle>0
\right\},
\qquad
B_r^-(x,\nu)
:=
\left\{
y\in B_r(x):
\langle\nu,y\rangle<0
\right\}.  
\end{equation}
\end{Definition}

The vectors $a,b\in V$ are called the \emph{(approximate) one-sided limits} of $u$ at $x$ with respect to $\nu$. 
Therefore,  for any $u\in BV^{\B}(\Omega)$, the triplet $(u^+,u^-,\nu_u): J_u\to V\times V\times S^{N-1}$ associated with the jump discontinuities on $J_u$  is well-defined, i.e.,
\[
x\in J_u
\quad\Longleftrightarrow\quad
\text{\eqref{eq:approxjumppoint} holds with }
(a,b,\nu)
=
(u^+(x),u^-(x),\nu_u(x)).
\]
Note that due to \eqref{BVBequiv} an analogous result hold for operators $\B$ of arbitrary order, replacing $u$ with the $(k-1)$-derivative $\nabla^{k-1}u$ of $u$ and considering $J_{\nabla^{k-1}u}$ (see \cite{ARS} for further details).  Following \cite{ARS}, for $u \in BV^{\mathcal B}(\Omega)$, we can decompose $\mathcal B u$ as follows 
\begin{align} \label{decompositionBu}
\mathcal B u = {\mathcal B}^a u + {\mathcal B}^s u = {\mathcal B}^a u + {\mathcal B}^s u \lfloor_{\Omega \setminus J_u} + {\mathcal B}^s u \lfloor_{J_u} =: {\mathcal B}^a u + {\mathcal B}^c u + {\mathcal B}^j u.
\end{align}

Note that, with an abuse of notation, we have defined  $\B^c u:=
{\mathcal B}^s u \lfloor_{\Omega \setminus J_u}$ which does not coincide with the classical notion of the Cantor part in $BV(\Omega)$. Indeed, in theory of $BV(\Omega)$, the Cantor part $D^{c}u$ is concentrated on the set of approximate continuity points of  $u$, while, for  $BV^{\B}(\Omega)$, with $\B$ being a $\mathbb C$-elliptic operator of arbitrary order $k$, it is not known whether
\begin{equation}\notag
|\B^c u|(S_u\setminus J_u)=0,
\end{equation}
where $S_u$ denotes the approximate discontinuity set of $u$ (for the definition of $S_u$ see, e.g., \cite[Definition 3.63]{AFP}). Therefore, the measure $\B^c u$ may, a priori, charge approximate discontinuity points that are not jump points.
In the sequel, we will use the following characterization of the jump part $\B^ju$ (see \cite[Corollary 2.8]{ARS} for operator of first order and \cite[Lemma 2.17]{ARS} for operators of arbitrary order)

\begin{lemma}
\label{lemma2.17ARS}
Let $\mathcal{B}$ be a $\mathbb{C}$-elliptic operator of order $k \geq 1$. If $u \in BV^{\mathcal{B}}(\Omega)$, then
\[
\mathcal{B}^j u = \mathbb{B}^k[\nu_{\nabla^{k-1}u}] \langle (\nabla^{k-1}u)^+ - (\nabla^{k-1}u)^-, \nu_{\nabla^{k-1}u}^{\otimes^{k-1}} \rangle_{k-1} \mathcal{H}^{n-1} \lfloor J_{\nabla^{k-1}u},
\]
where $\nu_{\nabla^{k-1}u}$ is an orientation of $J_{\nabla^{k-1}u}$, $\mathbb B^k$ is the symbol associated to $\B$ and defined in \eqref{mathbbAdef}. Here, we adopt the convention that for $k=1$, $\nabla^0u=u$.
\end{lemma}
Following the decomposition of $\B u$ given in \eqref{decompositionBu}, as in the classical theory of $BV$ spaces,
for a $\mathbb C$-elliptic operator $\B$, we can introduce the space $SBV^\B(\Omega)$ of special function of $\B$-bounded variations, consisting of functions $u\in BV^{\B}(\Omega)$ whose \textquotedblleft Cantor\textquotedblright\, part vanishes, i.e., 
    \begin{equation}
        \notag
        SBV^{\B}(\Omega)\coloneqq \{ u\in BV^{\B}(\Omega) \;:\; \B^{c} u = 0   \}.
    \end{equation}

\subsection{The Poincar\'{e} inequality}
We are going to prove a Poincar\'{e} inequality in the space $BV^\B(\Omega)$. To this end, a fundamental tool is \cite[Lemma 4.1.3]{Z} (see also \cite{T}) and for the reader's convience, we recall it below. 
\begin{lemma} {\rm{\cite[Lemma 4.1.3]{Z}}}
\label{lemma:Ziemerlemma4.1.3}
   Let $X_0$ be a normed linear space with norm $\|\cdot\|_0$ and let $X\subset X_0$ be a Banach space with norm $\|\cdot\|$. Assume that the norm $\|\cdot\|$ can be split as $\|\cdot\|_0+\|\cdot\|_1$, with $\|\cdot\|_1$ being a semi-norm in $X$. Moreover, assume that bounded sets in $X$ are precompact in $X_0$. Define $Y\coloneqq X\cap \{u\in X \;:\; \|x\|_1=0 \}$. If $L: X\to Y$ is a projection, there is a constant $C$ indipendent of $L$ such that 
        \begin{equation}
            \notag
            \|u-L(u)\|_0 \leq C\|L\|\|u\|_1, \;\;\; \hbox{for any } u\in X,
        \end{equation}
    with $\|L\|$ being the operator norm.
\end{lemma}
In our setting, we have that $X_0=L^1(\Omega)$ and $X= BV^{\B}(\Omega)$ endowed with the norm defined in \eqref{def:normBVB}, where $\Omega$ is a bounded open subset of $\R^N$ with Lipschitz boundary and $\B$ is an $\mathbb C$-elliptic homogeneous linear differential operator of order $k$ as in \eqref{Bdef}. The set $Y$ is then defined as
      \begin{equation}
          \notag
          Y\coloneqq BV^{\B}(\Omega)\cap \{u\in BV^{\B}(\Omega) : |\B u|=0     \}.
      \end{equation}
Note that $|\cdot|(\Omega)$ turns out to be a seminorm in $BV^{\B}(\Omega)$ and thanks to \cite[Corollary 2.16]{DEZ} bounded sets in $BV^{\B}(\Omega)$ are precompact in $L^1(\Omega)$. Thus, for any projection $L: BV^{\B}(\Omega)\to Y$, an application of Lemma \ref{lemma:Ziemerlemma4.1.3} yields to  
      \begin{equation}
            \notag
            \|u-L(u)\|_{L^1(\Omega)} \leq C\|L\||\B u|(\Omega), \;\;\; \hbox{for any } u\in BV^{\B}(\Omega).
        \end{equation}
In particular, for any bounded open convex $K$ containing the origin, we can introduce the projection operator  $L_K: BV^{\B}(K)\to Y$ defined as
    \begin{equation}
        \label{Lkproj}
        L_K(u) \coloneqq \frac{1}{{\mathcal H}^{n-1}(\partial K)} \int_{\partial K} u(x) d{\mathcal H}^{n-1}(x), \;\;\; u\in BV^{\B}(K).
    \end{equation}
Here, $u$ has to be understood in the sense of traces, which are defined \cite{DG24} for first-order operators and \cite{ARS} for $k$-th order operators. 

We now state a version of Lemma \ref{lemma:Ziemerlemma4.1.3}, which will be useful in the sequel.

\begin{Corollary}\label{cor:poincare}
Let $K$ be a bounded open convex set of $\R^N$. Let $L_k$ be the projection defined in \eqref{Lkproj} and let $u\in BV^{\B}(\rho K)$, with $\rho>0$. Then, there exist a positive constant $C$ only depending on $K$ and on $L_K$ such that 
      \begin{align}
      \label{poincareineq}
        \|u-L_{\rho K}(u)\|_{L^1(\rho K)} \leq  \rho^k C |\B u|(\rho K),
    \end{align}
    where $k$ is the order of the operator $\B$.
\end{Corollary}

\begin{proof}
 The inequality \eqref{poincareineq} follows from an homothety argument. Indeed, for $u\in BV^{\B}(\rho K)$, set $u_{\rho}(x)\coloneqq u(\rho x).$  By \eqref{Lkproj}, it immediately follows that  $L_{\rho K}(u) = L_{K}(u_\rho)$. This together with the classical change of variables implies that  $\|u-L_{\rho K}(u)\|_{L^1(\rho K)}= \rho^n \|u_\rho-L_{K}(u_\rho)\|_{L^1(K)}$. Moreover, we have that 
                    \begin{equation}
                        \label{eq:totvar+homarg}
                        |\B u_\rho|(\rho K) = \rho^{n+k}|\B u_{\rho}|(K).
                    \end{equation}
To prove \eqref{eq:totvar+homarg}, we use a density argument. Indeed, for $u\in BV^{\B}(\rho K)$, thanks to \cite[Lemma 6]{RS20}, there exists $(u_j)_j\in C^{\infty}(\rho K; V)\cap W^{\B, 1}(\rho K)$  such that $u_j$ $\B$-strictly converges to $u$, i.e., $u_j$ strongly converges to $u$ in $L^1(\rho K; V)$ and $|\B u_j|(\rho K)$ converges to $|\B u|(\rho K)$.  Then, by standard arguments and since $u_j$ is smooth, we deduce that
      \begin{align}
      \notag
          |\B u_j|(\rho K) &= \int_{\rho K}|\B u_j|(x)dx = \rho^n\int_{K}|\B u_j|(\rho x)dx = \rho^{n+k} \int_{K}|\B u_{j, \rho}|(x)dx = \rho^{n+k} |\B u_{j, \rho}u|(K), 
      \end{align}
with $u_{j, \rho}(x) \coloneqq u_j(\rho x)$. Now, passing to the limit as $j\to\infty$ and bearing in mind the $\B$-stricly convergences, we get \eqref{eq:totvar+homarg}. Therefore, an application of Lemma \ref{lemma:Ziemerlemma4.1.3} yields to, for $u\in BV^{\B}(\Omega)$,
    \begin{align}
        \|u-L_{\rho K}(u)\|_{L^1(\rho K)} &= \rho^n \|u_\rho-L_{K}(u_\rho)\|_{L^1(K)}\leq\rho^n C\|L_{K}\| |\B u_\rho|(K)=  \rho^{n+k} C\|L_{K}\| |\B u|(\rho K),\notag
    \end{align}
proving \eqref{poincareineq}, as desired.
\end{proof}

\section{The global method}\label{gm}

This section is devoted to the proof of a general integral representation theorem for  variational functionals defined in $BV^\B(\Omega)$. It is obtained as a generalization of the Global Method for Relaxation in $BV(\Omega)$ introduced in \cite{BFM}.

Let $\B$ be a linear, $\mathbb{C}$-elliptic partial differential operator of order $k$ as in \eqref{Bdef}. Let $\mathcal{F}: BV^{\B}(\Omega)\times {\mathcal A }(\Omega)\to [0, +\infty]$  be a functional satisfying the following assumptions
   \begin{itemize}
       \item[(H1)] ${\mathcal F}(u, \cdot)$ is the restriction of a Radon measure on $\mathcal A(\Omega)$ for any $u\in BV^{\B}(\Omega)$;
       \item[(H2)] ${\mathcal F}(\cdot, A)$ is $L^1(A; \R^d)$-lower semicontinuous for any $A\in {\mathcal A}(\Omega)$; 
       
       \item[(H3)] $\mathcal F(u; \cdot)$ is local, i.e. for all $A \in {\mathcal A}(\Omega)$,  if $v= u$ a.e. in $A$, then 
$$\mathcal F(v;A)= \mathcal F(u;A).$$
    \item[(H4)] there exists a positive constant $C>0$ such that for any $u\in BV^{\B}(\Omega)$ and $A\in {\mathcal A}(\Omega)$
                     \begin{equation}
                         \label{eq:growthcond}
                         \frac{1}{C}|\B u| (A) \leq {\mathcal F}(u; A) \leq C({\mathcal L}^N(A) + |\B u|(A)).
                     \end{equation}

    \end{itemize}
   In the sequel with an abuse of notation, we will identify $\mathcal F(u,\cdot)$ with its extension to the Borel subsets of $\Omega$, obtained via Theorem \ref{DGL}. In particular this extension inherits $(H3)$ for all Borel sets. This is, in fact, a consequence of the inner and outer regularity of the measures in the upper bound and in the lower bound. Moreover, by the monotonicity of the measure $\mathcal F(u;\cdot)$ and by the completeness of the measure $\mathcal L^N + |\mathcal B u|(\cdot)$, $\mathcal F(u;\cdot)$ is also a complete measure.
\color{black}
Our goal is to identify the bulk and the jump density of ${\mathcal F}$ satisfying Assumptions (H1)-(H4). To this end, we need to introduce a measure $m$ defined by 
          \begin{equation}
              \label{def:m}
              m(u; A) \coloneqq \inf \{ {\mathcal F}(v; A) \; :\; \text{supp}(u-v) \subset\subset A, v\in   BV^{\B}(\Omega) \}
          \end{equation}
We expect to prove the following integral representation
    \begin{align}
        \notag
        {\mathcal F}(u; A) &\coloneqq \int_A f(x, u, \nabla u, \dots, \nabla^k u)dx \\
        &\;\;\;\;\;\;\;\;\;\;\;\;\;\;\;\;\;\;+ \int_{J_{\nabla^{k-1}u}\cap A} h(x, u(x), \cdots, , \nabla^{k-2}u(x), (\nabla^{k-1}u)^{+},(\nabla^{k-1}u)^{-}, \nu_u) d{\mathcal H}^{N-1} \notag
    \end{align}
where 
  \begin{equation}
      \label{fdef}
      f(x_0, G_0, G_1, \dots, G_K) = \limsup_{\e\to 0} \frac{m(P^k_{x_0}(x-x_0), Q(x_0, \e))}{\e^N},
  \end{equation}
where $P^k_{x_0}$  is the polynomial defined in \eqref{Pkx0def},

\begin{align}\label{hdef}
h(x_0, G_0, G_1,\dots, G_{k-2}, c, d, \nu)=\limsup_{\varepsilon \to 0^+} \frac{m(P_{x_0}^{k-2}(x-x_0)+v_{J}(x), Q_\nu(x_0, \e))}{\varepsilon^{N-1}}
\end{align}
where $v_J$ is defined as in \eqref{v_J} 

First thing to prove is 
     \begin{equation}
         \notag
         \lim_{\e\to 0}\frac{{\mathcal F}(u; Q_{\nu}(x_0, \e))}{\mu(Q_{\nu}(x_0, \e)) } = \lim_{\e\to 0}\frac{m(u; Q_{\nu}(x_0, \e))}{\mu(Q_{\nu}(x_0, \e)) }
     \end{equation}
The following result will be used in the sequel. Its proof is omitted being very similar to the one in \cite[Lemma 4.1]{FHP}

\begin{lemma}
    Assume that Assumptions {\rm (H1)} and {\rm (H3)} hold, then for any $u\in $ and $A\in {\mathcal A}(\Omega)$
          \begin{equation}
              \notag
              \limsup_{\delta \to 0} m(u; A_\delta) \leq m(u; A),
          \end{equation}
    where $A_\delta \coloneqq \{x\in A\; : \; {\rm dist}(x, \partial A)>\delta\}$.
\end{lemma}

Set \begin{align}\label{mudef}
\mu\coloneqq {\mathcal L}^N + |\B u|.
\end{align}
Let $u\in BV^{\B}(\Omega)$. Fix $\delta>0$ and $A\in {\mathcal A}(\Omega)$, we define 
    \begin{align}
        \notag
        m^{\delta}(u; A)\coloneqq \inf\biggl\{ \sum_{i=1}^\infty m(u; B_i) \; : B_i\subset A\; \hbox{pairwise disjoint balls, diam}(B_i)\leq \delta, \mu\left(A\setminus \bigcup_{i=1}^\infty B_i\right)=0
        \biggr\}.
    \end{align}
Since $m^{\delta}(u; A)$ is increasing in $\delta$, the limit
          \begin{equation}
              \notag
              m^{\ast}(u; A) = \lim_{\delta\to 0}m^{\delta}(u; A) 
          \end{equation}
exists.

\begin{lemma}\label{Lemma 4.2}
    For all $A\in{\mathcal A}(\Omega)$, we have that 
          \begin{equation}
              \notag
              {\mathcal F}(u; A) = m^{\ast}(u; A).
          \end{equation}
\end{lemma}
\begin{proof}
Thanks to Assumption (H1) alongside with $m(u; A)\le {\mathcal F}(u; A)$, it follows immediately that $m^{\ast}(u; A)\le {\mathcal F}(u; A)$.\\
    We are going to prove the inequality $ {\mathcal F}(u; A) \leq  m^{\ast}(u; A).$  Fix $A \in \mathcal{A}(\Omega)$ and 
$\delta > 0$. Let $(B_i^{\delta})_i$ be the balls as in the definition of $m^{\delta}(u; A)$ such that
\begin{equation}
\label{eq:mlemdelta}
\sum_{i=1}^\infty m(u, B_i^{\delta}) \le m^{\delta}(u; A) + \delta.
\end{equation}
By defintion of $m$, there exists $v^\delta_i\in BV^{\B}(B_i^{\delta})$ such that $v_i^{\delta} = u$ in a neighborhood of $\partial B_i^{\delta}$ and 
     \begin{equation}
         \label{eq:Fleqm+delta}
         {\mathcal F}(v_i^{\delta}, B_i^{\delta}) \leq m(u, B_i^{\delta})+\delta. 
     \end{equation}
Set 
     \begin{equation}
         \notag
         v^{\delta}\coloneqq \sum_{i=1}^\infty v_i^\delta \chi_{B_i^\delta} + u \chi_{N_0^\delta},
     \end{equation}
where $N_0^\delta \coloneqq \Omega\setminus \bigcup_{i=1}^\infty B_i^\delta$. By construction, we have that $v^{\delta, n}$ is in $BV^{\B}(\Omega)$ and thanks to \eqref{eq:growthcond}, we also have that $v^{\delta}\in BV^{\B}(\Omega)$ and, thanks to \cite[Lemma 2.15 and 2.17]{ARS} (see also \cite[Corollary 4.21]{BDG}), we also have that  
         \begin{align}
            & \B v^\delta = \sum_{i=1}^\infty \B v_i^\delta\res {B_{i}^\delta} + \B u\res{N^\delta},\notag\\ 
            &|\B v^\delta|\res {N^\delta} = 0,\; \hbox{and}\; \mu(N^\delta)=0,\notag
         \end{align}
where $N^{\delta} \coloneqq A \cap N^\delta_0$. This, in turn, entails that
         \begin{equation}
             \notag
             {\mathcal F}(v^\delta; N^\delta) \leq C ({\mathcal L}^N(N^\delta) + |\B v^\delta|(N^\delta))=0.
         \end{equation}
         In fact, we are considering $\mathcal F$  already estended to the subsets of $\Omega$, is defined on $\mathcal A(\Omega)$
This combined with Assumption (H1),\eqref{eq:Fleqm+delta} and \eqref{eq:mlemdelta} implies that 
        \begin{align}
            {\mathcal F}(v^\delta; A) &= \sum_{i=1}^\infty {\mathcal F}(v^\delta_i; B_i^\delta) + {\mathcal F}(v^\delta; N^\delta)\notag\\
            &\leq \sum_{i=1}^\infty m(u; B_i^\delta)  +\delta \notag\\
            & \leq m^\delta (u; A) + 2\delta. \label{eq:Fleqmdelta}
         \end{align}
We claim that 
     \begin{equation}
         \label{eq:L1conv}
         v^\delta \to u\;\;\;\hbox{in}\;\;\; L^1(A).
     \end{equation}
Note that this allows us to conclude the proof. Indeed, thanks to Assumption (H2) as well as \eqref{eq:Fleqmdelta}, we have that 
      \begin{align}
          \notag
          {\mathcal F}(u; A)&\leq \liminf_{\delta\to 0} {\mathcal F}(v^\delta; A) \leq \liminf_{\delta\to 0} m^{\delta}(u; A) = m^{\ast} (u; A), 
      \end{align}
proving the proof of the inequality. \\
It remains now to prove \eqref{eq:L1conv}. The key tool is the Poincar\'{e} inequality.  Applying \eqref{poincareineq} and keeping in mind that 
$v_i^\delta = u$ in a neighborhood of $\partial B_i^\delta$, we have that 
         \begin{align}
             \notag
             \|v_i^\delta- u \|_{L^1(B_i^\delta)}=\|v_i^\delta- u \|_{L^1(B_i^\delta)}\leq C\;(\text{diam}(B_i^\delta))^k|\B (v_i^\delta- u)| (B_i^\delta).
         \end{align}   
Therefore, using the fact that ${\mathcal L}^N (A\setminus \cup B_i^\delta) \leq \mu((A\setminus \cup B_i^\delta))=0$ by definition of $m^\delta$
     \begin{align}
         \|v_i- u\|_{L^1(A)} &= \sum_{i=1}^\infty \|v_i^\delta- u\|_{L^1(B_i^\delta)}\notag\\
         & \leq C\;(\text{diam}(B_i^\delta))^k|\B (v_i^\delta- u)| \left(\bigcup_{i=1}^\infty B_i^\delta\right)\notag\\
         &\leq C \delta^k (|\B v^\delta|(A) + |\B u|(A)),\notag 
     \end{align}
which goes to $0$ as $\delta\to 0$ since $|\B v^\delta|(A)$ turns out to be bounded by Assumption (H3).
\end{proof}

The following result can be proven exactly as in \cite{BFM}. 
\medskip

\begin{Theorem}\label{asTheorem 4.3FHP}
Assume $(H1), (H2)$ and $(H3)$. Then for every $u\in BV^{\mathcal B}(\Omega)$ and
for all $\nu\in S^{N-1}$,
\[
\lim_{\varepsilon\to 0^+}
\frac{\mathcal F(u; Q_\nu(x_0,\varepsilon))}
     {\mu(Q_\nu(x_0,\varepsilon))}
=
\lim_{\varepsilon\to 0^+}
\frac{m(u; Q_\nu(x_0,\varepsilon))}
     {\mu(Q_\nu(x_0,\varepsilon))}
\]
for $\mu$–a.e.\ $x_0\in\Omega$, where 
$\mu := \mathcal{L}^N + |\mathcal B u|$.
\end{Theorem}

\begin{lemma}\label{Lemma 4.4FHP}
Assume $(H1)$-$(H4)$.  
Let $(v_\varepsilon)_\varepsilon\subset BV^\mathcal B(\Omega)$,
$u\in BV^{\mathcal B}(\Omega)$, $x_0\in\Omega$, $\nu\in S^{N-1}$,
and let $\lambda$ be a nonnegative Radon measure on $\Omega$.
Assume that there exists the following limit
\[
\lim_{\varepsilon\to 0^+}
\frac{ m(u;Q_\nu(x_0,\varepsilon)) }
     { \lambda(Q_\nu(x_0,\varepsilon)) }.
\]

Then
\[
\limsup_{\varepsilon\to 0^+}
\frac{ m(v_\varepsilon;Q_\nu(x_0,\varepsilon)) }
     { \lambda(Q_\nu(x_0,\varepsilon)) }
-\limsup_{\varepsilon\to 0^+}
\frac{ m(v_\varepsilon;Q_\nu(x_0,\varepsilon)) }
     { \lambda(Q_\nu(x_0,\varepsilon)) }\]
     \[
\leq\lim_{\delta \to 1^-}\lim_{\varepsilon\to 0^+}\frac{C}{\lambda(Q_\nu(x_0,\varepsilon))}\left\{\varepsilon^{N+1} +\mathcal |Bu|(Q_\nu(x_0,\varepsilon)\setminus    \overline{Q_\nu(x_0,\delta\varepsilon)})\right.
\]
\[
\qquad
+ |\mathcal B v_\varepsilon|(Q_\nu(x_0,\varepsilon)\setminus \overline{Q_\nu(x_0,\delta\varepsilon)})
+ \frac{1}{\varepsilon^k (1-\delta)^k}
 \int_{Q_\nu(x_0,\varepsilon)} |u - v_\varepsilon|\,dx
\]
\[
\left.\qquad
+ \dots+ \frac{1}{\varepsilon(1-\delta)}
 \int_{Q_\nu(x_0,\varepsilon)} |\nabla^{k-1} u - \nabla^{k-1} v_\varepsilon|\,dx
 \right\}.
\]
\end{lemma}

\begin{proof}
Fix $\delta\in(0,1)$ and let $\e>0$ be so small that $Q_{\nu}(x_0, \e)\subset\Omega$. Take $\phi\in C^\infty_{\rm c}(Q_{\nu}(x_0, \e))$ such that $\phi = 1$ in a neighborhood of $Q_{\nu}(x_0, \e)$ and 
     \begin{align}
         \notag
         \|\nabla ^m\phi \|_{\infty} \leq \frac{C}{{\e^{m+1}(1-\delta)^{m+1}}}, 
     \end{align}
for any $m\in\{0, \dots, k-1\}$, with $k$ being the order of the operator $\B$. Define
\[
w_\varepsilon =
\begin{cases}
\phi u + (1-\phi)v_\varepsilon &\text{in }Q_\nu(x_0,\varepsilon),\\
v_\varepsilon &\text{otherwise}.
\end{cases}
\]
Let $\tilde u$ be a function in the class of competitors of $m(u; Q_\nu(x_0,\varepsilon\delta))$ (see \eqref{def:m}) such that
\[
 {\mathcal F}(\tilde u;Q_\nu(x_0,\varepsilon\delta))\le \frac12 \varepsilon^{N+1} + m(u;Q_\nu(x_0,\varepsilon\delta)).
\]
By outer regularity (see (H1)), choose $\delta'\in (\delta, 1)$ with
\[
{\mathcal F}(\tilde u;Q_\nu(x_0,\varepsilon\delta'))
\le {\mathcal F}(\tilde u;Q_\nu(x_0,\varepsilon\delta))
+\frac12\varepsilon^{N+1}.
\]
Define
\[
\widetilde{w}_\varepsilon =
\begin{cases}
\tilde u &\text{in }Q_\nu(x_0,\varepsilon\delta),\\
w_\varepsilon &\text{on } \Omega\setminus Q_\nu(x_0,\varepsilon\delta).
\end{cases}
\]
Using locality (see (H3)), coercivity (see (H4)), we obtain 
    \begin{align}
        m(v_\e; Q_\nu(x_0, \e)) \leq m(\tilde u; Q_\nu(x_0, \e)) + \e^{N-1} +C\left( \e^N(1-\delta^N) + |\B w_\e|(Q_\nu(x_0, \e)\setminus\overline{Q_\nu(x_0, \e\delta)})   \right)  \notag
        \end{align}
Now, since 
     \begin{align}
         \B w_\e &= \B(\phi u) + \B((1-\phi)v_\e) = \phi \B u+ (1-\phi) \B v_\e + \sum_{\substack{|\beta|=k\\0\leq |\gamma|\leq k-1 }} c_{\beta\gamma} B_\beta \partial^{\beta-\gamma}\phi\partial^\gamma (u-v_\e)\notag\\
         &=\phi \B u+ (1-\phi) \B v_\e +[\symB, \phi](u-v_\e), \notag
     \end{align}
where $[\symB,\phi]$ denotes the communator which is a linear partial differential operator at most of order $k-1$ and whose coefficients depend only on $\|\nabla^m\phi\|_{\infty}$, for $m\in\{0, \dots, k-1\}$ and the principal symbol $\symB$. 
Then, 
   \begin{align}
       |\B w_\e|(Q_\nu(x_0, \e)\setminus\overline{Q_\nu(x_0, \e\delta)}) &\leq C\biggl( |\B u|(Q_\nu(x_0, \e)\setminus\overline{Q_\nu(x_0, \e\delta)})+|\B v_\e|(Q_\nu(x_0, \e)\setminus\overline{Q_\nu(x_0, \e\delta)})  \notag\\
       &\quad    + \frac{1}{\e^k(1-\delta)^k} \int_{Q_\nu(x_0, \e)}  |u-v_\e|dx +\cdots+ \frac{1}{\e(1-\delta)} \int_{Q_\nu(x_0, \e)}  |\nabla^{k-1} u-\nabla^{k-1} v_\e|dx  \biggr). \notag
   \end{align}
\end{proof}

\begin{Corollary}
\label{Cor45FHP}
Assume ${\rm (H1)}$-${\rm (H4)}$.  
Let $v,u\in BV^\B(\Omega)$, $x_0\in\Omega$, 
$\nu\in S^{N-1}$, and $\lambda$ a nonnegative Radon measure.
Assume that there exist the following limits
\[
\lim_{\varepsilon\to 0^+}\frac{m(u;Q_\nu(x_0,\varepsilon))}
{\lambda(Q_\nu(x_0,\varepsilon))},
\qquad
\lim_{\varepsilon\to 0^+}\frac{m(v;Q_\nu(x_0,\varepsilon))}
{\lambda(Q_\nu(x_0,\varepsilon))}.
\]
Then,
\[
\left|\limsup_{\varepsilon\to 0^+}
\frac{ m(v;Q_\nu(x_0,\varepsilon)) }
     { \lambda(Q_\nu(x_0,\varepsilon)) }
-\limsup_{\varepsilon\to 0^+}
\frac{ m(u;Q_\nu(x_0,\varepsilon)) }
     { \lambda(Q_\nu(x_0,\varepsilon)) }\right|\]
     \[\quad\leq
\lim_{\delta \to 1^-}\lim_{\varepsilon\to 0^+}\frac{C}{\lambda(Q_\nu(x_0,\varepsilon))}\left(\varepsilon^{N+1} + \mathcal |Bv|(Q_\nu(x_0,\varepsilon)\setminus    \overline{Q_\nu(x_0,\varepsilon\delta)})\right.
\]
\[
\left.\qquad
+ |\mathcal B u|(Q_\nu(x_0,\varepsilon)\setminus \overline{Q_\nu(x_0,\varepsilon\delta)})
+ \frac{1}{\varepsilon^k (1-\delta)^k}
 \int_{Q_\nu(x_0,\varepsilon)} |u - v|\,dx  + \cdots + \frac{1}{\varepsilon (1-\delta)}
 \int_{Q_\nu(x_0,\varepsilon)} |\nabla^{k-1} u - \nabla^{k-1} v|\,dx  
 \right).
\]
\end{Corollary}

This theorem gives an integral representation in $SBV^\B$.

\begin{Theorem}
\label{thm:globalmethod}
Assume ${\rm (H1)}$-${\rm (H4)}$.  
For every $u\in SBV^\B(\Omega)$ and every $A\in\mathcal{A}(\Omega)$,
\[
{\mathcal F}(u;A)
=
\int_A f(x,u,\nabla u,\dots, \nabla^k u)\,dx
\;+\;
\int_{J_{\nabla^{k-1} u}\cap A}
h(x,u, \dots, \nabla^{k-2}u, (\nabla^{k-1}  u)^+,(\nabla^{k-1}  u)^-,\nu_{\nabla^{k-1} u})\,d\mathcal{H}^{N-1},
\]
where $f$ and $h$ are given in \eqref{fdef} and \eqref{hdef}, respectively.

\end{Theorem}

\medskip

\begin{proof}
\color{black}
	We begin by proving that, for $\mathcal L^N$- a.e. $x_0 \in \Omega$,
	\begin{align}\label{fproof}
	\frac{d \mathcal F(u;\cdot)}{d \mathcal L^N} (x_0)=	
	f(x_0, u(x_0), \nabla u(x_0), \dots, \nabla^k u(x_0)).
	\end{align}
Let $P^k_{x_0}(x-x_0)$ be as in \eqref{Pkx0def} and let $x_0$ be a fixed point in $\Omega$ satisfying the following properties
\begin{align}
&\lim_{\varepsilon \to 0^+}\frac{1}{\varepsilon^k}\ave_{Q(x_0,\varepsilon)}|u(x) - P^k_{x_0}(x-x_0)| \, d x = 0;\label{aproxdif}\\
&\lim_{\varepsilon \to 0^+}\frac{1}{\varepsilon^{k-m}}\ave_{Q(x_0,\varepsilon)}|\nabla^m u(x) - \nabla^m P^k_{x_0}(x-x_0)| \, d x = 0, \; \hbox{for } m=1,\dots, k-1;\label{aproxdifnabla}\\
& \mathcal B^a u(x_0):=\lim_{\varepsilon \to 0^+}\frac{1}{\varepsilon^N}\B u(Q(x_0,\varepsilon)) = B^k(\nabla^ku(x_0)), \quad \lim_{\varepsilon \to 0^+}\frac{1}{\varepsilon^N}|\B^su|(Q(x_0,\varepsilon)) = 0;\label{DuDsu}\\
& \lim_{\varepsilon \to 0^+}\frac{1}{\varepsilon^N}\int_{Q(x_0,\varepsilon)}|B^k (\nabla^k u(x))-B^k(\nabla^ku(x_0))|dx =0;  \label{Lebpoint}      \\
& \frac{d \mathcal F(u;\cdot)}{d \mathcal L^N} (x_0) =
\lim_{\varepsilon \to 0^+}\frac{\mathcal F(u;Q(x_0,\varepsilon))}{\varepsilon^N} =
\lim_{\varepsilon \to 0^+}\frac{m(u;Q(x_0,\varepsilon))}{\varepsilon^N};\label{Fmu}\\
& \frac{d \mathcal F(P^k_{x_0}(x-x_0);\cdot)}{d \mathcal L^N} (x_0) =
\lim_{\varepsilon \to 0^+}\frac{m(P^k_{x_0}(x-x_0);Q(x_0,\varepsilon))}{\varepsilon^N}.\label{Fmva}
\end{align}

Note that such a point $x_0$ esists. Indeed, the approximate differentiability property \eqref{aproxdif} follows from \cite[Lemma 2.3 and Lemma 2.4]{R} together with the embedding theorems  \cite[Theorem 3.7]{DG24}. We emphasize that $\nabla^k u(x_0)$ denotes the coefficient of the highest-order term in the approximate  polynomial $P^k_{x_0}$, and should not, in general, be interpreted as a Sobolev derivative. Observe that through the jet map $B^k$, it determines the density of the absolutely continuous part of $\mathcal Bu$.  The property \eqref{aproxdifnabla} follows from the embedding theorem \cite[Theorem 3.7]{DG24} combined with standard differentiability results in Sobolev spaces, \cite[eq. (2.3)]{R}, \cite[Lemma 2.3]{R} which relies on Calderon-Zygmund' s theorem 11.1 in \cite{CZ}, and also \cite[the proof of Lemma 6.1]{R} (see \cite[Section 6]{EG} for the first-order case and \cite[Theorem 2.1]{FHP} for the second-order case). 

Furthermore, the identities \eqref{DuDsu} follow from \cite[Lemma 2.4]{R}, while \eqref{Lebpoint} is a consequence of both \cite[Lemma 2.4]{R} and \cite[Corollary 2.23]{AFP}. Finally, in establishing \eqref{Fmu} and \eqref{Fmva}, we also make use of Theorem \ref{asTheorem 4.3FHP}. All the above properties hold for $\mathcal L^N$-almost every $x_0\in\Omega$. 

Fix $x_0\in \Omega$ as above. Let $\delta \in (0,1)$ and let $\varepsilon >0$ be small enough so that $Q(x_0,\varepsilon)\subset \Omega$. Given the definition of the density $f$ in \eqref{fdef}, due to \eqref{Fmu} and \eqref{Fmva}, we want to show that

\begin{align}\label{toestimate}
	\lim_{\varepsilon \to 0^+}
	\frac{m(P^k_{x_0}(x-x_0);Q (x_0, \varepsilon))}
	{\mathcal L^N (Q (x_0, \varepsilon))}- \lim_{\varepsilon \to 0^+}\frac{m(u; Q (x_0, \varepsilon))}{\mathcal L^N (Q (x_0, \varepsilon))} = 0.
\end{align}

By Lemma \ref{Lemma 4.4FHP}, we know that 
   \begin{align}\lim_{\varepsilon \to 0^+}
	&\frac{m(P^k_{x_0}(x-x_0);Q (x_0, \varepsilon))}
	{\mathcal L^N (Q (x_0, \varepsilon))}- \lim_{\varepsilon \to 0^+}\frac{m(u; Q (x_0, \varepsilon))}{\mathcal L^N (Q (x_0, \varepsilon))}\notag\\
      &\quad \leq\lim_{\delta \to 1^-}\lim_{\varepsilon\to 0^+}\frac{C}{\e^N}\left\{\varepsilon^{N+1} +\mathcal |\mathcal B u|(Q_\nu(x_0,\varepsilon)\setminus    \overline{Q_\nu(x_0,\delta\varepsilon)})\right. \notag\\
      &\qquad
+ |B^k (\nabla^k u(x_0))|(Q_\nu(x_0,\varepsilon)\setminus \overline{Q_\nu(x_0,\delta\varepsilon)})
+ \frac{1}{\varepsilon^k (1-\delta)^k}
 \int_{Q_\nu(x_0,\varepsilon)} |u - P^k_{x_0}(x-x_0)|\,dx\notag\\
 &\left.\qquad
+ \dots+ \frac{1}{\varepsilon(1-\delta)}
 \int_{Q_\nu(x_0,\varepsilon)} |\nabla^{k-1} u - \nabla^{k-1} P^k_{x_0}(x-x_0)|\,dx
 \right\}. \notag
   \end{align}
In view of Lemma \ref{FM}, we have that 
    \begin{align}\notag
       0&\leq \lim_{\delta \to 1^-}\lim_{\varepsilon\to 0^+}\mathcal |\mathcal Bu|(Q_\nu(x_0,\varepsilon)\setminus   \overline{Q_\nu(x_0,\delta\varepsilon)}) \leq\lim_{\delta \to 1^-} |B^k(\nabla^ku(x_0))|(1-\delta^N) = 0.
    \end{align}
This combining with \eqref{aproxdif} and \eqref{aproxdifnabla} yields to \eqref{toestimate}.\\

Next, we want to prove that, for $\mathcal H^{N-1}$-a.e $x_0 \in J_{\nabla^{k-1}u}$,
\begin{align*}
\frac{d \mathcal F(u; \cdot)}{d \mathcal H^{N-1}\lfloor{J_{\nabla^{k-1}u}}}(x_0) = h(x_0, u(x_0), \nabla u(x_0), \dots, \nabla^{k-1}u(x_0),  (\nabla^{k-1}u)^+(x_0), (\nabla^{k-1}u)^{-}(x_0), \nu_{\nabla^{k-1}u(x_0)}).
\end{align*}
For simplicity of notation we denote by $\nu$ the unit vector $\nu_{\nabla^{k-1}u(x_0)}$ and by $v_J$ the function defined in \eqref{v_J} with $(c,d)=(\nabla^{k-1}u^+(x_0),\nabla^{k-1}u^-(x_0))$, i.e. 
$$v_J(x) = v_{\nabla^{k-1}u^+(x_0),\nabla^{k-1}u^-(x_0),\nu}(x)=\frac{1}{(k-1)!}\left\{
\begin{array}{ll}
\langle \nabla^{k-1}u^+(x_0),(x-x_0)^{\otimes^{k-1}}\rangle &\hbox{ if } (x-x_0)\cdot \nu> 0,\\
\langle \nabla^{k-1}u^-(x_0), (x-x_0)^{\otimes^{k-1}}\rangle &\hbox{ if }(x-x_0)\cdot \nu \leq 0. 
\end{array}\right.
$$
For $\mathcal H^{N-1}$-almost every $x_0\in J_{\nabla^{k-1}u}$, and for $k$ sufficiently large, the following properties hold:

\begin{align}
&\quad \lim_{\varepsilon \to 0^+}\frac{1}{\varepsilon^{k-1}}\ave_{Q_{\nu}(x_0,\varepsilon)}|u(x) -P^{k-2}_{x_0}(x-x_0)- v_J(x)| \, dx =0;\label{jumppt} \\
&\quad \lim_{\varepsilon \to 0^+}\frac{1}{\varepsilon^{k-2}}\ave_{Q_{\nu}(x_0,\varepsilon)}|\nabla u(x) -\nabla P^{k-2}_{x_0}(x-x_0)- \nabla v_J(x)| \, dx =0\notag\\
&\quad \lim_{\varepsilon \to 0^+}\frac{1}{\varepsilon^{k-3}}\ave_{Q_{\nu}(x_0,\varepsilon)}|\nabla^2 u(x) -\nabla^2 P^{k-2}_{x_0}(x-x_0)- \nabla^2 v_J(x)| \, dx =0\notag\\
&\quad \;\;\;\ \vdots\notag\\
&\quad \lim_{\varepsilon \to 0^+}\frac{1}{\varepsilon}\ave_{Q_{\nu}(x_0,\varepsilon)}|\nabla^{k-2} u(x) -\nabla^{k-2} P^{k-2}_{x_0}(x-x_0)- \nabla^{k-2} v_J(x)| \, dx =0\label{jumpptk-2}\\
 &\quad \lim_{\varepsilon \to 0^+}\ave_{Q_{\nu}(x_0,\varepsilon)}|\nabla^{k-1}u(x) -\nabla^{k-1}v_J(x)| \, dx =0\label{jumpptk-1}\\
&\quad \lim_{\varepsilon \to 0^+}\frac{1}{\varepsilon^{N-1}}|\mathcal B u|(Q_{\nu}(x_0,\varepsilon)) =\lim_{\varepsilon \to 0^+}\frac{1}{\varepsilon^{N-1}}|B^k({D^k u})|(Q_{\nu}(x_0,\varepsilon)) = |{\mathbb B}[\nu](\nabla^{k-1}u^{+}(x_0) -\nabla^{k-1}u^{-}(x_0))|;\label{Du}\\
&\quad \frac{d \mathcal F(u;\cdot)}{d \mathcal H^{N-1}\lfloor{J_{\nabla^{k-1}u}}} (x_0) =
\lim_{\varepsilon \to 0^+}\frac{\mathcal F(u;Q_{\nu}(x_0,\varepsilon))}{\varepsilon^{N-1}} =
\lim_{\varepsilon \to 0^+}\frac{m(u;Q_{\nu}(x_0,\varepsilon))}{\varepsilon^{N-1}};\label{jFmu}\\
&\quad \frac{d \mathcal F(v_J;\cdot)}{d \mathcal H^{N-1}\lfloor{J_{\nabla^{k-1}u}}} (x_0) =
\lim_{\varepsilon \to 0^+}\frac{m(v_J;Q_{\nu}(x_0,\varepsilon))}{\varepsilon^{N-1}};\label{Fmvj}.\\
\label{Uip}
\end{align}

Note that such a point  $x_0\in \Omega$ exists. Indeed, estimates from \eqref{jumppt} to \eqref{jumpptk-1} can be deduced from \cite[Theorem 2.4, Remark 2.3, Corollary 2.8, and Lemma 2.17]{ARS}, together with \cite{D}. On the other hand, for the readers' convenience we deduce \eqref{jumppt} in Corollary \ref{cor:jumppt} as a consequence of Proposition \ref{prop:estAFP}, the subsequent properties for the derivatives, until \eqref{jumpptk-2} follow analogously. Finally, \eqref{Du} follows from Lemma \ref{lemma2.17ARS}.
\smallskip

Therefore, fix $x_0\in \Omega$ satisfying  properties \eqref{jumppt}-\eqref{Fmvj}. Applying Corollary \ref{Cor45FHP} with $v=v_j$ and $\lambda= H^{N-1}\lfloor{S_{\nabla^{k-1}u}}$, we claim that
\begin{align}
   &\left|\limsup_{\varepsilon\to 0^+}\frac{ m(P^{k-2}_{x_0}(\cdot-x_0)+ v_J;Q_\nu(x_0,\varepsilon)) }
     { \e^{N-1} }-\limsup_{\varepsilon\to 0^+}
\frac{ m(u;Q_\nu(x_0,\varepsilon)) }
     { \e^{N-1} }\right| \notag\\
     &\leq \lim_{\delta \to 1^-}\lim_{\varepsilon\to 0^+}\frac{C}{\e^{N-1}}\Big(\varepsilon^{N+1} + \mathcal |Bv_j|(Q_\nu(x_0,\varepsilon)\setminus    \overline{Q_\nu(x_0,\varepsilon\delta)})\notag\\
    & \qquad
+ |\mathcal B u|(Q_\nu(x_0,\varepsilon)\setminus \overline{Q_\nu(x_0,\varepsilon\delta)})
+ \frac{1}{\varepsilon^k (1-\delta)^k}
 \int_{Q_\nu(x_0,\varepsilon)} |u(x) - P^{k-2}_{x_0}(x-x_0)- v_J(x)|\,dx\notag\\
 &\qquad+ \frac{1}{\varepsilon^{k-1} (1-\delta)^k}
 \int_{Q_\nu(x_0,\varepsilon)} |\nabla u(x) - \nabla P^{k-2}_{x_0}(x-x_0)- \nabla v_J(x)|\,dx  + \cdots \color{black}\notag\\
 &\qquad + \frac{1}{\varepsilon (1-\delta)}
 \int_{Q_\nu(x_0,\varepsilon)} |\nabla^{k-1} u - \nabla^{k-1} v_J|\,dx  
 \Big)=0.\label{estimate1} 
\end{align}
Indeed, making use of \eqref{jumppt}-\eqref{Du}, we deduce that 
    \begin{align}\notag
        0\leq \lim_{\delta \to 1^-}\lim_{\varepsilon\to 0^+}\frac{ |\mathcal B u|(Q_\nu(x_0,\varepsilon)\setminus \overline{Q_\nu(x_0,\varepsilon\delta)})}{\e^{N-1}}\leq \lim_{\delta \to 1^-} \|\mathbb B\| |[\nabla^{k-1}u](x_0)|(1-\delta^{N-1}) = 0.
    \end{align}
This combining with \eqref{jFmu} and \eqref{Fmvj} proves \eqref{estimate1}, concluding the proof of the theorem.

\end{proof}

We conclude this section, giving the sketch of Theorem \ref{thm:GMSDB}.

\begin{proof}[Sketch of Theorem \ref{thm:GMSDB}]
 The proof follows the same arguments of those in Theorem \ref{thm:globalmethod}, with minor modifications which we outline below.  Since the set where $\mathcal F$ is defined is $SBV^\B \times L^1$, the class of competitors (see \eqref{compSDB}) used in definition \eqref{m4.3} is different from the class in \eqref{mudef}, but the proof of the lemmata proceeds along the same lines. 
 In the same spirit, for the identification of the bulk density $f$, we further require $x_0\in \Omega$ satisfying \eqref{aproxdif}-\eqref{Fmva} to be a Lebesgue point of $U$, i.e.,
     \begin{equation}
         \notag
         \lim_{\e\to 0}\frac{1}{\e^N}\int_{Q(x_0, \e)} |U(x)- U(x_0)|dx = 0. 
     \end{equation}
As for the surface density $g$,  we additionally require $x\in J_{\nabla^{k-1}u}$ satisfying \eqref{jumppt}-\eqref{Uip} to be such that 
         \begin{equation}
             \notag
              \lim_{\e\to 0}\frac{1}{\e^{N-1}}\int_{Q(x_0, \e)} |U(x)|dx = 0.
         \end{equation}
   Under these additional conditions, the identification of $f$ and $g$ proceeds exactly as in the proof of Theorem \ref{thm:globalmethod}.        
\end{proof}

\begin{Remark}\label{rem:caseBVxL1}
An analogous result holds in the space $BV^{\B}(\Omega)\times L^1(\Omega; ImB^k)$, and its proof follows the same lines as that of Theorem~\ref{thm:GMSDB}.  In this setting, the functional $\mathcal F$ as in \eqref{def:introF} is defined on $BV^{\B}(\Omega)\times L^1(\Omega;\operatorname{Im} B^k)$, and the class of admissible competitors for the minimization problem $m$ in \eqref{m4.3} is defined on the same space. Therefore, the bulk and surface densities $f$ and $h$ appearing in the corresponding integral representation are those associated with this minimization problem.
\end{Remark}

\section{Applications}
\label{sec:applications}

This section provides some applications of Theorem \ref{thm:GMSDB}, regarding the integral representations of some special class of functionals. Let $\Omega$ be an open and bounded domain of $\R^N$. 

Let $\B$ be a $k$-th order $\mathbb{C}$-elliptic differential operator defined in \eqref{Bdef}.
Throughout this section, we assume that $V$ is isomorphic to $\R^d$ and $W$ is a subspace of $\mathbb R^{d} \otimes E_k(\R^N)$. We first restrict our attention to the class $SBV^{\B}(\Omega)$ defined by 
     \begin{equation}\notag
         SBV^{\B}(\Omega)\coloneqq \{ u\in BV^{\B}(\Omega) \; : \;  \B^{c} u = 0   \},
     \end{equation}
with $\B^c u$ being the Cantor contribution of the measure $\B u\in {\mathcal M}(\Omega; W)$. We consider the following functional

 \begin{equation}
     \label{eq:functionalF0}
     F_0(u; A) \coloneqq
     \begin{cases}
         \int_A f_0(x, u, \nabla u, \dots, \nabla^{k-1}u, \B u)dx\\
            \;\;\;\;\;\;\;\; + \int_{J_{\nabla^{k-1}u}\cap A} g_0(x, u, \nabla u, \dots, \nabla^{k-2}u, (\nabla^{k-1}u)^+, (\nabla^{k-1}u)^-, \nu_{\nabla^{k-1} u}) d{\mathcal H}^{N-1}, & u\in SBV^{\B}(\Omega);\\
            \\
            +\infty, & \hbox{otherwise},
     \end{cases}
 \end{equation}
 where the densities $f_0$ and $g_0$ are such that
     \begin{itemize}
         \item[(G1)] $f_0: \Omega \times \mathbb R^d \times \mathbb R^{d\times N}\times \mathbb R^{d} \otimes E_2(\R^N)\times \dots \times W \to [0, \infty)$ is measurable in $x$ and continuous in all other variables and there exists a constant $C>0$ such that
              \begin{equation} \label{eq:crescitaf_0}
                  \frac{1}{C} |\Lambda|\leq f_0(x, a, \xi,\dots, \xi^{k-1}, \Lambda ) \leq C(1+|\Lambda|),
              \end{equation}
        for all $(x, a, \xi,\dots, \xi^{k-1}, \Lambda ) \in\Omega \times \mathbb R^d\times \mathbb R^{d\times N}\times \mathbb R^{d\times N^2}\times \dots \times W$;
        \item[(G2)] $g_0: \Omega \times \mathbb R^d \times \mathbb R^{d\times N}\times \mathbb R^d \otimes E_2(\mathbb R^{N})\times \dots \times \mathbb R^{d}\otimes E_{k-1}(\mathbb R^N)\times \mathbb R^{d}\otimes E_{k-1}(\mathbb R^N)\times \mathbb S^{N-1} \to [0, \infty)$ is continuous and there exists a constant $C>0$ such that
           \begin{equation}\label{eq:crescitag0}
               \frac{1}{C} |\xi^{k-1}-\eta^{k-1}| \leq g_0(x, a, \xi, \dots,\xi^{k-2}, \xi^{k-1}, \eta^{k-1}, \nu) \leq C(1+ |\xi^{k-1}-\eta^{k-1}|),
           \end{equation}
         for all $(x, a, \xi,\dots, \xi^{k-1}, \eta^{k-1},\nu) \in\Omega \times \mathbb R^d \times \mathbb R^{d\times N}\times \mathbb R^{d}\otimes E_2(\mathbb R^N)\times \dots \mathbb R^{d}\otimes E_{k-1}(\mathbb R^N) \times \mathbb R^{d}\otimes E_{k-1}(\mathbb R^N)\times \mathbb S^{N-1}$. Moreover, $g_\e$ is symmetric, i.e., 
               \begin{equation}
                   \notag
                   g_0(x, a, \xi, \dots,\xi^{k-2}, \xi^{k-1}, \eta^{k-1}, \nu) = g_0(x, a, \xi, \dots,\xi^{k-2}, \eta^{k-1}, \xi^{k-1}, -\nu). 
               \end{equation}
     \end{itemize}
    
     \color{black}
     The functional ${\mathcal F}: SBV^{\B}(\Omega)\times {\mathcal A}(\Omega)\to [0, \infty)$ is defined by 
           \begin{align}
               \label{relaxeddef}
               {\mathcal F}(u; A) \coloneqq \inf \biggl\{ \liminf_{n\to\infty} F_0(u_n; A): \;\; u_n\to u \; \hbox{in } W^{k-1, 1}(A; \R^d),\;\; \B u_n {\mathcal L}^N\wsto \B u \; \hbox{in } {\mathcal M}(A; W) \biggr\}.
           \end{align}
The aim is to show that the functional ${\mathcal F}$ defined above satisfies all the assumptions of Theorem \ref{thm:globalmethod}.

\subsection{Integral representation in the \texorpdfstring{$SBV^{\B}$}{SBH} setting}
\begin{lemma}\label{lemma:lsc}
Let $\Omega$ be a bounded open subset in $\mathbb R^N$. Let $F_0$ be as in \eqref{eq:functionalF0} with $f_0$ and $g_0$ as in \eqref{eq:crescitaf_0} and \eqref{eq:crescitag0}, respectively, and let  $\mathcal F$ be as in \eqref{relaxeddef}.
    For every $u\in SBV^{\B}(\Omega)$, $A\in{\mathcal A}(\Omega)$ and every sequence $(u_n)_n \subset SBV^{\B}(\Omega)$ such that $u_n\to u$ in $W^{k-1, 1}(\Omega;\mathbb R^d)$ and $\B u_n {\mathcal L}^N\wsto \B u$ in ${\mathcal M}(\Omega; W)$,
         \begin{equation}
             \notag
             {\mathcal F}(u; A)\leq \liminf_{n\to\infty} {\mathcal F}(u_n; A).
         \end{equation}
\end{lemma}
\begin{proof}
Fix a sequence $(u_n)_n\subset SBV^{\B}(A)$ as in the statement
For every $(u_n)_n$, 
then there exists a subsequencea
$(u_{n,k})_k\subset SBV^{\B}(A)$
such that $u_{n,k}\to u_n$ in $W^{k-1,1}(A; \R^d)$ and
$\B u_{n,k}{\mathcal L}^N\stackrel{*}{\rightharpoonup} \B u_n$ in ${\mathcal M}(A; W)$ as $k\to\infty$ and
\[
\liminf_{k\to\infty} F_0(u_{n,k};A)
\leq
\mathcal{F}(u_n;A)+\frac{1}{n}.
\]
By diagonalizing, we find sequences
$v_n:=u_{n,k_n}$ such that $v_n\to u$ in $W^{k-1, 1}(A;\mathbb R^d)$,
$\B v_n{\mathcal L}^N\stackrel{*}{\rightharpoonup} \B u$ in ${\mathcal M}(A;W)$ as $n\to\infty$, and
\[
\liminf_{n\to\infty}
F_0(v_n;A)
\leq
\liminf_{n\to\infty}
\mathcal{F}(u_n;A)
\]
and thus
\[
\mathcal{F}(u;A)
\leq
\liminf_{n\to\infty}
\mathcal{F}(u_n;A),
\]
concluding the proof.
\end{proof}

From the definition of $\mathcal F$ in \eqref{relaxeddef} it immediately follows that
\begin{lemma}\label{lemma:loc}
The functional $\mathcal{F}$ is local, i.e.,
for all $A \in \mathcal{A}(\Omega)$, if $u=v$ $\mathcal{L}^N$-a.e.\ $x\in A$,
then
\[
\mathcal{F}(u;A)=\mathcal{F}(v;A).
\]
\end{lemma}

\begin{lemma}\label{lemma:growth}
Let $f_0$ and $g_0$ be two densities satisfying Assumption {\rm (G1)} and {\rm (G2)}.
For every $u\in SBV^{B}(\Omega)$ and
for every $A\in\mathcal{A}(\Omega)$, we have
\[
\frac{1}{C}
|\B u|(A)
\leq \mathcal{F}(u;A)
\leq C\left(
\mathcal{L}^N(A)
+|\B u|(A)
\right),
\]
with $C>0$. Moreover, for every
$u\in BV^{\B}(\Omega)$ the functional
$\mathcal{F}(u;\cdot)$ is the restriction to
$\mathcal{A}(\Omega)$ of a Radon measure.
\end{lemma}

\begin{proof}
We note that Assumptions (G1) and (G2) imply that
\[
\frac{1}{C}|\B u|(A)
\leq F_0(u;A)
\leq C\left(
\mathcal{L}^N(A)+|\B u|(A)
\right)
\]
for every $u\in SBV^{\B}(\Omega)$
and $A\in\mathcal{A}(\Omega)$.
For any $u\in SBV^{\B}(\Omega)$ and any
$\delta>0$, we can find
$u_n\in SBV^{\B}(\Omega)$ such that
$u_n\to u$ in $W^{k-1, 1}(\Omega;\mathbb R^d)$,
$\B u_n\stackrel{*}{\rightharpoonup}\B u$ in ${\mathcal M}(\Omega;W)$,
and
\[
\mathcal{F}(u;A)
\geq \liminf_{n\to\infty}
F_0(u_n;A)-\delta.
\]

On one hand, this implies that
\[
\mathcal{F}(u;A)
\geq \liminf_{n\to\infty}
\frac{1}{C}|\B u_n|(A)-\delta
\geq \frac{1}{C}|\B u|(A)-\delta,
\]
and letting $\delta\to0$ we have
\begin{equation}
\label{eq:5.3}
\mathcal{F}(u;A)
\geq \frac{1}{C}|\B u|(A).
\end{equation}
To prove the upper bound, we consider the sequence
$(u_n)_n$ constructed in \cite[Proposition 2.13]{DEZ} which, in particular, satisfies
$u_n \to u$ in $W^{k-1,1}(A;\R^d)$,
$|\B u_n|(A) \to |\B u|(A)$. 
Then, we have
\[
\begin{aligned}
\mathcal{F}(u;A) \leq \liminf_{n\to\infty}
F_0(u_n;A)\leq \liminf_{n\to\infty}
C\left(
\mathcal{L}^N(A)+|\B u_n|(A)
\right)\leq C\left(
\mathcal{L}^N(A)+|\B u|(A)
\right).
\end{aligned}
\]

Finally, we will prove that for
$u\in SBV^{\B}(\Omega)$,
$\mathcal{F}(u;\cdot)$ is the restriction
to $\mathcal{A}(\Omega)$ of a Radon measure.
We will apply a variant of the De Giorgi-Letta criterion as stated in \cite[Lemma 2.3]{FHP} and recalled in the Appendix for the readers' convenience as Lemma \ref{DGL}. 
Since item (ii) follows directly from the fact that
$F_0(u,\cdot)$ is a Radon measure and
item (iii) follows from the growth condition that
we have just proved, it only remains to prove
that for any open sets
$A,B,C\in\mathcal{A}(\Omega)$ with
$\overline{A}\subset B\subset C$ we have
\[
\mathcal{F}(u;C)
\leq
\mathcal{F}(u;C\setminus\overline{A})
+\mathcal{F}(u;B).
\]

To see this, for $\varepsilon>0$ we choose
$v_n\in BV^{\B}(C \setminus \bar A)$ and
$w_n\in BV^{\B}(B)$ as in the
definition of $\mathcal{F}(u;\cdot)$
(perhaps along a subsequence) so that
\begin{equation}\label{eq: 5.4}
\lim_{n\to\infty}
F_0(v_n,C\setminus\overline{A})
\leq
\mathcal{F}(u;C\setminus\overline{A})
-\varepsilon.
\end{equation}
and
\[
\lim_{n\to\infty}
F_0(w_n,B)
\leq
\mathcal{F}(u;B)-\varepsilon.
\]

We will use a slicing argument in order to construct
(up to a possible extraction of a subsequence) a sequence
$(u_n)_n\subset BV^{\B}(C)$
as in the definition of $\mathcal{F}(u;\cdot)$
so that
\[
\liminf_{n\to\infty} F_0(u_n;C)
\leq
\lim_{n\to\infty}
F_0(v_n,C\setminus\overline{A})
+
\lim_{n\to\infty}
F_0(w_n,B).
\]

Let $\delta>0$ be so small that 
\[
S_{\delta}\coloneqq \{x\in B: \hbox{dist}(x, A)<\delta\}\subset \subset B.
\]

Given $k\in\mathbb{N}$, we can decompose
$S_\delta\setminus A$ into a disjoint union of strips, i.e.
\[
S_\delta\setminus A
=\bigcup_{i=1}^{\ell} L_{i,\ell},
\]
where
\[
L_{i,\ell}
=\left\{
x\in S_\delta :
\frac{(i-1)\delta}{\ell}
<\operatorname{dist}(x,A)
\leq\frac{i\delta}{\ell}
\right\}.
\]

By coercivity of $F_0$, we have
\[
\sup_n |\B v_n|(C\setminus\overline{A})
+\sup_n |\B w_n|(B)
\leq M
\]
for some $M<\infty$, and thus
\[
\sup_n \sum_{i=1}^{\ell}
\left(
|\B v_n|+|\B w_n|
\right)(L_{i,\ell})
\leq M.
\]

We remark that since there are only finitely many
values of $i$ and infinitely many values of $n$,
there must be some fixed $i$ such that
\[
\left(
|\B v_n|+|\B w_n|
\right)(L_{i,\ell})
\leq \frac{M}{\ell}
\]
for infinitely many $n\in\mathbb{N}$.
Thus, for any $\ell$, there are
an $i_\ell\in\{1,\ldots,\ell\}$ and a subsequence
$(n_j^{(\ell)})\subset(n)$ such that
\[
\left(
|\B v_{n}|+|\B w_{n}|
\right)(L_{i_\ell,\ell})
\leq \frac{M}{\ell}, 
\]
for infinitely many $n\in\mathbb{N}$.
Thus, for any $\ell$, there are
an $i_\ell\in\{1,\ldots,\ell\}$ and a subsequence
$(n_j^{(\ell)})\subset(n)$ such that
\[
\left(
|\B v_{n_j^{(\ell)}}|
+
|\B w_{n_j^{(\ell)}}|
\right)(L_{i_\ell,\ell})
\leq \frac{M}{\ell},
\quad \forall\,j,\ell\in\mathbb{N}.
\]

We consider a smooth cutoff function
$\phi_\ell\in C_c^\infty(B;[0,1])$ such that
$\{0<\phi_\ell<1\}\subset L_{i_k,k}$,
$\phi_\ell(x)=0$ if
$\operatorname{dist}(x,A)
\leq \frac{i_\ell-1}{\ell}\delta$,
$\phi_\ell(x)=1$ if
$\operatorname{dist}(x,A)
\geq \frac{i_\ell}{\ell}\delta$ and
\[
\|\nabla\phi_\ell\|_\infty\leq C\ell,\cdots, \|\nabla^{k-1}\phi_\ell\|_\infty\leq C\ell^{k-1},
\|\nabla^{k}\phi_\ell\|_\infty\leq C\ell^k.
\]

For $x\in C$, we define
\[
u_{j,\ell}
=
\phi_\ell v_{n_j^{(\ell)}}
+
(1-\phi_\ell)w_{n_j^{(\ell)}}.
\]

Then, we have
\[
\begin{aligned}
F_0(u_{j,\ell};C)
&\leq
F_0(v_{n_j^{(\ell)}};C\setminus\overline{A})+
F_0(w_{n_j^{(\ell)}};B)
+
F_0(u_{j,\ell};L_{i_\ell,\ell}),
\end{aligned}
\]

and the last term is bounded by
\[
\begin{aligned}
F_0(u_{j,\ell};L_{i_\ell,\ell})
&\leq C\Bigg(
\mathcal{L}^N(L_{i_\ell,\ell})
+\ell^k\int_{L_{i_\ell,\ell}}
|v_{n_j^{(\ell)}}-w_{n_j^{(\ell)}}|\,dx +\ell^{k-1}\int_{L_{i_\ell,\ell}}
|\nabla v_{n_j^{(\ell)}}-\nabla w_{n_j^{(\ell)}}|\,dx+\cdots\\
&\qquad
+\left(|\B v_{n_j^{(\ell)}}|
+|\B w_{n_j^{(\ell)}}|\right)(L_{i_\ell,\ell})
\Bigg)\\
&\leq C\Bigg(
\frac{1}{\ell}
+\ell^{k}\int_{L_{i_\ell,\ell}}
|v_{n_j^{(\ell)}}-w_{n_j^{(\ell)}}|\,dx +\ell^{k-1}\int_{L_{i_\ell,\ell}}
|\nabla v_{n_j^{(\ell)}}-\nabla w_{n_j^{(\ell)}}|\,dx + \cdots \notag\\
&\qquad + \ell\int_{L_{i_\ell,\ell}}
|\nabla^{k-1} v_{n_j^{(\ell)}}-\nabla^{k-1} w_{n_j^{(\ell)}}|\,dx
\Bigg).
\end{aligned}
\]

Since $v_n\to u$ and $w_n\to u$ in
$W^{k-1,1}(B\setminus\overline{A})$, for any $\ell$
we can choose an element $n_{j_\ell}^{(\ell)}$ of
$n_j^{(\ell)}$ so that the map
$\ell\mapsto n_{j_\ell}^{(\ell)}$ is increasing and
\[
\int_{B\setminus\overline{A}}
\left|
v_{n_{j_\ell}^{(\ell)}}-w_{n_{j_\ell}^{(\ell)}}
\right|\,dx
=o(1/\ell^k) , \cdots, \int_{B\setminus\overline{A}}
\left|
\nabla^{k-1} v_{n_{j_\ell}^{(\ell)}}
-\nabla^{k-1} w_{n_{j_\ell}^{(\ell)}}
\right|\,dx
=o(1/\ell).
\]
With this choice, we have that
\[
\lim_{\ell\to\infty}
F_0(u_{j_\ell,\ell};L_{i_\ell,\ell})=0.
\]

Since $v_n\to u$ in
$W^{k-1,1}(C\setminus\overline{A})$,
$w_n\to u$ in $W^{k-1,1}(B)$ and
$\B v_n\stackrel{*}{\rightharpoonup}\B u$
in $C\setminus\overline{A}$,
$\B w_n\stackrel{*}{\rightharpoonup}\B u$
in $B$, we must have that
$u_{j_\ell,\ell}\to u$ in $W^{k-1,1}(C)$ and
$\B u_{j_\ell,\ell}\stackrel{*}{\rightharpoonup} \B u$
in $C$. Thus, by definition of
$\mathcal{F}(u;\cdot)$, we conclude
\[
\begin{aligned}
\mathcal{F}(u;C)
&\leq \liminf_{\ell\to\infty}
F_0(u_{j_\ell,\ell};C)\\
&\leq \lim_{\ell\to\infty}
F_0(v_{n_{j_\ell}^{(\ell)}};
C\setminus\overline{A})\\
&\quad+
\lim_{\ell\to\infty}
F_0(w_{n_{j_\ell}^{(\ell)}};B)\\
&\leq
\mathcal{F}(u;C\setminus\overline{A})
+\mathcal{F}(u;B)-2\varepsilon.
\end{aligned}
\]
Sending $\varepsilon\to0$, we are done.
\end{proof}

Consequently since $\mathcal F$ in \eqref{relaxeddef} satisfies the assumptions of the Global Method, as proved in  Lemmata \ref{lemma:growth}, \ref{lemma:loc} and \ref{lemma:lsc}, Theorem \ref{thm:globalmethod} applies to $\mathcal F$ and the following integral representation on $SBV^{\B}(\Omega)$ holds.

\begin{Theorem}[Integral representation in $SBV^{\B}$]\label{thm:intrepSBVB}
Let $\B$ be a linear, $\mathbb{C}$-elliptic partial differential operator of order $k$. Let $F_0$ be the functional defined by \eqref{eq:functionalF0} satisfying Assumptions {\rm (G1)}--{\rm (G2)} and let $\mathcal F$ its relaxed functional defined by \eqref{relaxeddef}. Then, for every $u \in SBV^{\B}(\Omega)$ and every $A \in \mathcal{A}(\Omega)$, $\mathcal F$ admits the following representation
\begin{align}
  \mathcal{F}(u;A) &= \int_A f(x,u,\nabla u,\dots,\nabla^k u)\,dx\notag\\
&\qquad+ \int_{J_{\nabla^{k-1} u} \cap A} h(x,u,\nabla u, \dots, \nabla^{k-2}u, (\nabla^{k-1} u)^+, (\nabla^{k-1} u)^-, \nu_{\nabla^{k-1}u})\, d\mathcal{H}^{N-1},  \notag
\end{align}
where the bulk density $f$ and the interfacial density $h$ are given respectively by \eqref{fdef} and \eqref{hdef}.
\end{Theorem}

We can push further the investigation of the integral representation of the relaxed functional $\mathcal F$ of $F_0$ under stronger assumptions of the densities $f_0$ and $g_0$, in the same spirit of \cite[Theorem 4.2.2]{BFM}.\\
Let $f_0: \Omega \times\mathbb R^d \times \mathbb R^{d\times N}\times \mathbb R^{d} \otimes E_2(\R^N)\times \dots\times \mathbb R^{d} \otimes E_{k-1}(\R^N) \times W \to [0, \infty)$ be a continuous function satisfying \eqref{eq:crescitaf_0}, and
   \begin{itemize}
       \item[(I1)] for any $\e>0$ there exists $\delta>0$ such that for any $x,y \in \Omega$ and $a,b \in \mathbb R^d$,  $\xi,\eta \in \mathbb R^{d\times  N},\xi^2, \eta^2\in \mathbb R^{d} \otimes E_2(\R^N), \cdots, \xi^{k-1},\eta^{k-1}\in \mathbb R^{d} \otimes E_{k-1}(\R^N)$,  with $|x-y|+|a-b|+\dots |\xi^{k-1}-\eta^{k-1}|<\delta$ we have
            \begin{equation}
                \notag
                |f_0(x, a, \xi, \dots, \xi^{k-1}, \Lambda)-f_0(y, b, \eta, \dots, \eta^{k-1}, \Lambda)| \leq C\varepsilon (1+|\Lambda|),
            \end{equation}
             for every $\Lambda  \in W$;
            \end{itemize}
    Let $g_0:\Omega \times\mathbb R^d \times \mathbb R^{d\times N}\times \mathbb R^{d}\otimes E_2(\mathbb R^N) \times \dots \times \mathbb R^{d}\otimes E_{k-1}(\mathbb R^N) \times\mathbb R^{d}\otimes E_{k-1}(\mathbb R^N) \times \mathbb S^{N-1} \to [0, \infty) $ be a continuous function satisfying \eqref{eq:crescitag0} and
    \begin{itemize}
    \item[ (I2)] for any $\e>0$ there exists $\lambda>0$ such that for any $x,y \in \Omega$ and $a,b \in \mathbb R^d$,  $\xi,\zeta \in \mathbb R^{d\times  N},\cdots, \xi^{k-2},\zeta^{k-2}\in \mathbb R^{d}\otimes E_{k-2}(\mathbb R^N)$,  with $|x-y|+|a-b|+\dots |\xi^{k-2}-\zeta^{k-2}| + |z|<\delta$ we have
         \begin{equation}
             \notag
             |g_0(x, a, \xi, \dots, \xi^{k-1}-z, \eta^{k-1}-z, \nu)-g_0(y, b, \zeta, \dots,\xi^{k-1}, \eta^{k-1}, \nu)| \leq C\varepsilon (1+|\xi^{k-1}-\eta^{k-1}|),
         \end{equation}
    for any $\xi^{k-1}, \eta^{k-1}\in \mathbb R^{d} \otimes E_{k-1}(\R^N)$.
   \end{itemize}

   \begin{Theorem}\label{repcont}
Let $\Omega$ be a bounded open set. Let $f_0$ and $g_0$ be continuous functions as above, satisfying \eqref{eq:crescitaf_0}, $(I1)$ and \eqref{eq:crescitag0}, $(I2)$, respectively.  Let $F_0:SBV^\B(\Omega)\to[0,+\infty)$ be as in \eqref{eq:functionalF0} and let $\mathcal{F}: SBV^{\B}(\Omega)\times\mathcal{A}(\Omega)\to [0, +\infty)$ be as in \eqref{relaxeddef}, then
\begin{align}\label{Frepcont}
\mathcal F(u,\Omega)=& \int_\Omega \widetilde{f}(x, u(x), \nabla u(x), \dots, \nabla^ku(x)) dx \notag\\
&\quad+ \int_{J_{\nabla^{k-1}u}\cap\Omega} \widetilde{h}(x, u(x),\nabla u(x) \dots, (\nabla^{k-1}u)^{+}(x), (\nabla^{k-1}u)^{-}(x), \nu_{\nabla^{k-1}u}]) dx,
\end{align}
where, for all $x_0\in \Omega$, $a\in \R^d$, $\xi\in\R^{d\times N}, \xi^2\in \mathbb R^{d} \otimes E_2(\R^N), \dots, \xi^{k-2}\in \mathbb R^{d} \otimes E_{k-2}(\R^N)$, $\xi^{k-1}, \eta^{k-1}\in \mathbb R^{d} \otimes E_{k-1}(\R^N)$, $\nu\in \mathbb S^{N-1}$ and $\xi^k\in \mathbb R^d \otimes E_k(\mathbb R^N)$
    \begin{align}
        \notag
        &\widetilde{f}(x_0, a, \xi, \cdots, \xi^{k-1}, \xi^k) \coloneqq \limsup_{\e\to 0}\inf\biggl\{\int_{Q} f_0(x_0, a, \xi, \cdots, \xi^{k-1}, \B v(x))dx\notag\\
        &\quad+\frac{1}{\e}\int_{J_{\nabla^{k-1} u}\cap Q} g_0(x_0, a, \xi, \cdots, \xi^{k-2}, \xi^{k-1}+\e (\nabla^{k-1} v)^{+}(x), \xi^{k-1}+\e (\nabla^{k-1} v)^{-}(x), \nu_v(x))d{\mathcal H}^{N-1}(x)\notag\\
        &\quad\quad\quad: v\in SBV^{\B}(Q)\;\;\hbox{and }  \hbox{supp}(v-\widetilde{P})\subset\subset Q 
        \biggr\}\label{def:tildef}
    \end{align}
and 
    \begin{align}
        \notag
        &\widetilde{h}(x_0, a, \xi, \cdots,\xi^{k-2}, \xi^{k-1}, \eta^{k-1}, \nu) \coloneqq \limsup_{\e\to 0}\inf\biggl\{\int_{Q_\nu} \e f_0\left(x_0, a, \xi, \cdots, \xi^{k-2}, \nabla^{k-1} v(x), \frac{\B v(x)}{\e}\right)dx\notag\\
        &\quad+\int_{J_{\nabla^{k-1} v}\cap Q_\nu}     g_0(x_0, a, \xi, \cdots,\xi^{k-2},(\nabla^{k-1} v)^{+}(x), (\nabla^{k-1} v)^{-}(x), \nu_v(x))  d{\mathcal H}^{N-1}(x)\notag\\
        &\quad\quad\quad: v\in SBV^{\B}(Q)\;\;\hbox{and }  \hbox{supp}(v-\widetilde{v}_J)\subset\subset  Q_\nu
        \biggr\}\label{def:tildeh}
    \end{align}
where $\widetilde{P}(x)\coloneqq \frac{1}{k!}<\xi^k, x^{\otimes k}>$ and  $\widetilde{v}_J$ is given by 
                  \begin{align}
                      \notag
                      \widetilde{v}_J (x)\coloneqq \frac{1}{(k-1)!}
                      \begin{dcases}
                          \langle \xi^{k-1}, x^{\otimes k}\rangle, &   \hbox{if } x\cdot \nu >0,\\
                          \langle \eta^{k-1}, x^{\otimes k}\rangle, &   \hbox{if } x\cdot \nu \leq0.
                      \end{dcases}
                  \end{align}
    \end{Theorem}

The proof of Theorem \ref{repcont} follows the same arguments of \cite[Theorem 4.2.2]{BFM}, thus here we just outlined the main steps. To this end, we need to introduce a new minimum problem. For $u\in SBV^{\B}(\Omega)$ and $A\in {\mathcal A}_\infty(\Omega)$, define
  \begin{align}\notag
      m_0(u; A) \coloneqq \inf \biggl\{ F_0(v;A):
      v\in SBV^{\B}(A), \;\;\; \hbox{supp}(v-u)\subset\subset A
      \biggr\}.
  \end{align}

\begin{lemma}\label{lem:m=m_0}
    Let $f_0$ and $g_0$ be continuous functions as above, satisfying \eqref{eq:crescitaf_0}, $(I1)$ and \eqref{eq:crescitag0}, $(I2)$, respectively.  Let $F_0:SBV^\B(\Omega)\to[0,+\infty)$ be as in \eqref{eq:functionalF0} and let $\mathcal{F}: SBV^{\B}(\Omega)\times\mathcal{A}(\Omega)\to [0, +\infty)$ be as in \eqref{relaxeddef}. Then, for any $u\in SBV^{\B}(\Omega)$ and $A\in{\mathcal A}_\infty(\Omega)$, we have that
         \begin{equation}
             \notag
             m_0(u; A)= m(u; A), 
         \end{equation}
         where $m$ is the minimum problem in \eqref{def:m} associated to $\mathcal F$ in \eqref{relaxeddef}.
\end{lemma}
The proof of Lemma \ref{lem:m=m_0} is entirely similar to the one of \cite[Lemma 4.2.2]{BFM}, hence it is omitted. However, it is worth to observe that the proof  hinges on Lemma \ref{lem:DeGiorgislicing} as well as the following lemma which allows us to manipulate the boundary values and whoose proof is very similar to the one of \cite[Lemma A.2]{DEZ}.
\begin{lemma}\label{lem:DeGiorgislicing} 
Let \(F_0\) and \(\mathcal F : BV^\B(\Omega)\times\mathcal A(\Omega) \to [0,+\infty] \) be the functionals as in \eqref{eq:functionalF0} and \eqref{relaxeddef}.
Let \(u\in SBV^\B(\Omega)\) and let \((v_n)_n\subset SBV^\B(\Omega)\) be such that \( v_n \wsto u \hbox{ in }BV^\B(\Omega). \) Then, for every \(A\in\mathcal A_\infty(\Omega)\), there exists a sequence \((w_n)_n\subset BV^\B(\Omega)\) such that \( w_n\wsto u \hbox{ in }BV^\B(\Omega), \) \( w_n=u \hbox{ in a neighborhood of }\partial A, \) and \[ \limsup_{n\to\infty} \mathcal F(w_n;A) \le \liminf_{n\to\infty} \mathcal F(v_n;A). \] 
\end{lemma}
    
    \begin{proof}[Proof of Theorem \ref{repcont}]
    The assumptions on $f_0$ and $g_0$ ensure that $\mathcal F$ admits the representation provided by Theorem \ref{thm:intrepSBVB}, with the integrands $f$ and $h$ being defined as in \eqref{fdef} and \eqref{hdef} associated to $\mathcal F$ in \eqref{relaxeddef},
respectively. To obtain more explicit formulas for these integrands, we apply Lemma \ref{lem:m=m_0}, which allows us to replace $m$ with $m_0$ in the definitions \eqref{fdef} and \eqref{hdef}. Indeed, replacing $m$ by $m_0$ in  \eqref{fdef} along with performing the standard change of variables  $y=\tfrac{x-x_0}{\varepsilon}$ and introducing the rescaled test functions defined by
        \begin{equation}
            \notag
            \widetilde{v}(y)\coloneqq \frac{1}{\e^{k}}\left(v(x_0+\e y) - P^{k-1}_{0}(\e y)   \right),
        \end{equation}
    where $P_{0}^{k-1}$ is the polynomial defined in \eqref{Pkx0def} with $x_0=0$, we get \eqref{def:tildef}. Similarly, to obtain the corresponding formula for $h$, we perform the same change of variables and introduce the rescaled test functions defined by
     \begin{equation}
         \notag
         \widetilde{v}(y)\coloneqq \frac{1}{\e^{k-1}}\left(v(x_0+\e y) - P^{k-2}_{0}(\e y)   \right),
     \end{equation}
     with $P^{k-2}_{0}$ being defined as above.
\end{proof}

\subsection{Integral representation in $BV^{\B}$}
Theorem \ref{thm:globalmethod} provides an integral representation result in $SBV^{\B}(\Omega)$. If, in addition, we suppose that the abstract functional $\mathcal F$ satisfying Assumptions (H1)--(H4) is invariant under ${k-1}$-th polynomial translations of $u$ and continuous with respect to the area strictly convergence, we also have a representation result on the entire space $BV^{\B}(\Omega)$, as stated in the next result.

\begin{Corollary} \label{cor:affineinva}
 Let $\B$ be a linear, $\mathbb C$-elliptic partial differential operator of order $k$. Let $\Omega$ be a bounded domain in $\mathbb{R}^N$ and let $A\in{\mathcal A}_{\infty}(\Omega)$.  Let $\mathcal F$ be a functional satisying Assumptions {\rm (H1)}--{\rm (H4)}.
 Assume also that for any $u\in BV^{\B}(\Omega)$ and $R\in \hbox{ker}\B$, we have that 
    \begin{equation}
        \label{invariance}
        \mathcal F(u + R, A) = \mathcal F(u;A).
    \end{equation}
\color{black}
\begin{itemize}
    \item[{\rm (i)}] Then, for any $u\in SBV^{\B}(\Omega)$ and $A\in{\mathcal A}(\Omega)$, we have, with a slight abuse of notation, 
\[
{\mathcal F}(u;A) = \int_A f\left(x, \frac{d\B u}{d{\mathcal L}^N}(x)\right)\,dx
+ \int_{J_{\nabla^{k-1} u}\cap A} h(x, (\nabla^{k-1}u)^{+},(\nabla^{k-1}u)^{-} ,\nu_{J_{\nabla^{k-1}u}})\, d\mathcal{H}^{N-1},
\]
where $f$ in \eqref{fdef} is constant with respect to $\mathbb R^d \otimes E_0(\mathbb R^N) \times \dots \mathbb R^d \otimes E_{k-1}(\mathbb R^N)$, and $h$ in \eqref{hdef} is constant with respect $E_0(\mathbb R^N) \times \dots \mathbb R^d \otimes E_{k-2}(\mathbb R^N)$. In particular, if $u \in W^{\B, 1}(\Omega)$, the functional reduces to
\begin{equation}\label{repSob}
{\mathcal F}(u;A) = \int_A f(x,\B u(x))\,dx.
\end{equation}

\item[{\rm (ii)}] If $u \in BV^\B(\Omega)$, and the energy density
$f:\Omega \times \mathbb R^d \times E_k(\mathbb R^N)\to [0,+\infty)$, given by \eqref{repSob}, satisfies
\begin{itemize}
\item[$(f_1)$] $f$ is Lipschitz continuous in the second variable uniformly in $x \in \overline \Omega$;
\item[$(f_2)$] there exists a modulus of continuity $\omega:[0,+\infty)\to [0,+\infty)$ such that
$|f(x, \Lambda)- f(y,\Lambda)|\leq \omega (|x-y|)(1+ |\Lambda|)$
for all $x,y\in \overline{\Omega}$ and $\Lambda\in \mathbb R^d \otimes E_k(\mathbb R^N).$
\end{itemize}
Then, it results that
\[
{\mathcal F}(u;A) \leq \int_A f\left(x,\frac{d\B u}{d{\mathcal L}^N}(x)\right)\,dx + \int_A f^\infty\left(x, \frac{d \B u}{d \B^s u}(x)\right) d |\B^s u|(x).
\]
where $f^\infty(x,\Lambda):=\lim_{t \to +\infty}\frac{f(x, t \Lambda)}{t}$ for $x \in \overline{\Omega}$ and $\Lambda \in \mathbb R^d \otimes E_k(\mathbb R^N)$.

\item[{\rm (iii)}] If $\Omega$ has also a Lipschitz boundary  and if $\mathcal F$  is continuous  with respect to the area-strict convergence in \eqref{def:areaconv} relative to $\mu_n = \mathcal B u_n$ and $\mu= \mathcal B u$. Then, for every $u \in BV^{\B}(\Omega)$ 
\[
\mathcal F(u) = \int_\Omega f\left(x,\frac{d\B u}{d{\mathcal L}^N}(x)\right)\,dx
+ \int_\Omega f^\infty\!\left(x, \frac{d(\B u)^s}{d|(\B u)^s|}(x)\right) d|(\B u)^s|(x).
\]
\end{itemize}

\end{Corollary}

\begin{proof}[Sketch of the proof]
The first item is a consequence of Theorem \ref{thm:globalmethod}.
For what concerns the integral representation in $BV^\B(\Omega)$, due to the fact that $f$ satisfies $(f_1)$ and $(f_2)$ and since $\mathcal F$ is lower semicontinuous in $BV^\B(\Omega)$ it results that $f=\mathcal Q_{B}f$ in view of \cite[Theorem 1.3]{DEZ} and the integral representation on $W^{1,\B}(\Omega)$. Hence it suffices to invoke again \cite[Theorem 1.3]{DEZ} to get the desired inequality on $BV^{\mathcal B}(\Omega)$.

For what concerns the last item, we observe that
\eqref{repSob} provides a representation of $\mathcal F$ on
$W^{\B,1}(\Omega)$.
For any $u\in BV^\B(\Omega)$, we can use \cite[Proposition 2.13]{DEZ}, to construct a sequence
$(u_n)_n\subset W^{\B,1}(\Omega)$
so that $u_n\to u$ in $L^1$ and
$\B u_n\to \B u$ area-strictly.
Therefore, by the assumption on the area-strict continuity of $\mathcal F$, we have
\begin{equation}
\mathcal F(u)
=\lim_{n\to\infty} \mathcal F(u_n).
\label{eq:6.5}
\end{equation}

On the other hand, the functional
\[
\begin{aligned}
u\in BV^{\B}(\Omega)
\longmapsto\;& I(u):=
\int_\Omega f\left(x,\frac{d\B u}{d{\mathcal L}^N}(x)\right)\,dx +\int_\Omega
f^\infty\left(
x,\frac{d(\B u)^s}{d|(\B u)^s|}(x)
\right)\,d|(\B u)^s|(x)
\end{aligned}
\]
is area-strict continuous on $BV^{\B}(\Omega)$ \cite[Theorem 3.2 - Proof of Upper Bound]{DEZ} and agrees
with $\mathcal F$ on $W^{\B,1}(\Omega)$, therefore
\[
\lim_{n\to\infty}\mathcal F(u_n)
=
\lim_{n\to\infty}I(u_n)=I(u).
\]

This, along with \eqref{eq:6.5}, yields to 
\[
\begin{aligned}
\mathcal{F}(u;\Omega)
={}&\int_\Omega f\left(x,\frac{d\B u}{d{\mathcal L}^N}(x)\right)\,dx+\int_\Omega
f^\infty\left(
x,\frac{d(\B u)^s}{d|(\B u)^s|}(x)
\right)\,d|(\B u)^s|(x)
\end{aligned}
\]
for every $u\in BV^{\B}(\Omega)$.
\end{proof}

\begin{Remark}
Assumption \eqref{invariance} generalizes the corresponding condition introduced in \cite{CVG}.  
\end{Remark}

\subsection{Integral representation in the space $SD_\B$}\label{subsec:SD}

In this section, we establish an integral-representation result in the space of structured deformations  $SD_{\B}$ (see Definition \ref{def:SD} for its definition). Throughout, we restrict our analysis to differential operators $\B$ of order $k=1$ or $k=2$. This restriction is motivated by the approximation results for $BV^{\mathcal B}$ fields that play a crucial role in proving the integral representation. In the first-order case, the relevant approximation result is due to Alberti \cite{alberti}. For second-order operators, an analogous result was first established in $\R^2$ in \cite{FLP} and then extended to bounded domains $\Omega\subset\R^N$ in arbitrary dimension in \cite{FHP}. \\
We start proving an approximation result for operator $\B$ of first order. To this end, we recall the approximation result of Alberti \cite{alberti} in the version quoted in  \cite[Theorem~2.8]{CF1997}.

\begin{Theorem}[{Alberti's Theorem}]\label{Al}
Let $G \in L^1(\Omega; \mathbb R^{d{\times} N})$. Then, there exist a function $f \in SBV(\Omega; \mathbb R^d)$, a Borel function $\beta\colon\Omega\to\mathbb R^{d{\times} N}$, and a constant $C_N>0$ depending only on $N$ such that
\begin{equation*}
Df = G \,{\cL}^N + \beta {\mathcal H}^{N-1}\res J_f, \qquad
\int_{\Omega\cap J_f}  |\beta| \, d {\mathcal H}^{N-1} \leq C_N \lVert G\rVert_{L^1(\Omega)}.
\end{equation*}
\end{Theorem}

Moreover, we will use also the following results.

\begin{lemma}[{\cite[Lemma~2.9]{CF1997}}]\label{ctap}
Let $u \in BV(\Omega; \mathbb R^d)$. Then there exist piecewise constant functions $\bar u_n\in SBV(\Omega;\mathbb R^d)$  such that $\bar u_n \to u$ in $L^1(\Omega; \mathbb R^d)$ and
\begin{equation*} 
|Du|(\Omega) = \lim_{n\to \infty}| D\bar u_n|(\Omega) = \lim_{n\to \infty} \int_{\Omega\cap J_{\bar u_n}} |[\bar u_n]|\; d{\mathcal H}^{N-1}.
\end{equation*}
\end{lemma}

\begin{Theorem}[{\cite[Theorem~1.2]{Silhavy}}]\label{silhavy_app}
 Let  $(g,G)\in BV(\Omega;\mathbb R^N)\times L^1(\Omega;\mathbb R^{d\times N})$. Then,  there exists a sequence $(g_n)_n \subset SBV(\Omega;\mathbb R^N)$ such that
$g_n \wsto g$ in $BV(\Omega;\mathbb R^N)$ and $\nabla g_n =G$.
\end{Theorem}
Then, we can prove the approximation result in $BV^\B$, with $\B$ being an operator of order $k=1$.

\begin{Theorem}\label{ourappthm}
 Let $\mathcal{B}$ be a $\mathbb C$-elliptic linear homogeneous differential operator of first order, whose symbol in jet form is $B$.  Let ${\rm Im}B\subseteq W$ be as in \eqref{ImBk}. Let  $(g,G)\in BV^{\mathcal B}(\Omega; V)\times L^1(\Omega;{\rm Im} B)$.
 Then,  there exists a sequence $(g_n)_n \subset SBV^{\B}(\Omega;V)$ such that
$g_n \wsto g$ in $BV^{\mathcal B}(\Omega;V)$ and $\mathcal B^a g_n \wsto G$ in $\mathcal M(\Omega;{\rm Im}B)$.
\end{Theorem}

\begin{proof}[Proof] The proof follows the same arguments as in the $SBV$ case, by applying Alberti's Theorem and a piecewise constant approximation of $L^1$ fields, but taking into account the presence of $\mathcal B$. Namely let $B$ the unique symbol associated to $\mathcal B$ for the jet form $\mathcal B g$ as in \cite[(2.8)]{ARS}, cf \eqref{jetformBu}.  Consider the symbol $B^{\dagger}$ defined as the Moore-Penrose pseudoinverse of $B$. 

Recalling \eqref{ImBk}, ${\rm Im} B$ is a subset of $W$.
For every $G \in L^1(\Omega;{\rm Im}B)$, consider the jet form $B^{\dagger} G$ which is now an element of $V\otimes E_1(\mathbb R^N)= V\otimes \mathbb R^N$, where $B^{\dagger} $ is the Moore-Penrose generalized inverse. Recall that for any matrix $B\in \R^{m\times N}$ there is a unique matrix $B^{\dagger} \in \R^{N\times m}$ such that 
     \begin{equation}
         \notag
         B^{\dagger} BB^{\dagger} = B^{\dagger}, \;\; BB^{\dagger}B=B, \;\; (BB^{\dagger})^* = BB^{\dagger} \;\; (B^{\dagger}B)^{*}= B^{\dagger} B.
     \end{equation}
 For the properties of $B^{\dagger}$ see, e.g.,  \cite{CM}. Then Theorem \ref{Al} allows to find a field $f \in SBV(\Omega; V)$  such that $\nabla f =B^{\dagger}G$. Then, again applying the jet form to both elements, it results $\mathcal B^a f = B \nabla f = B B^\dagger G = G $ since $G \in {\rm Im }B$.

In particular, since $\mathcal B^a g= B\nabla g \in L^1(\Omega;{\rm Im}B)$  we can define $u_n= g + h - h_n$,
where by Alberti's Theorem $h \in SBV(\Omega; V)$ (with $B\nabla h(x) \in {\rm Im}B \subseteq W$ for a.e. $x \in \Omega$),  is such that
\begin{align}\label{817bis} 
B \nabla h =  \mathcal{B}^ah = G - \mathcal{B}^ag, \quad \quad |D h|(\Omega) \leq C \Vert B^\dagger G \Vert_{L^1(\Omega)}  + C \Vert B^\dagger B\nabla g \Vert_{L^1(\Omega)} \leq C'( \Vert G \Vert_{L^1(\Omega)}  +\Vert\nabla g \Vert_{L^1(\Omega)}). 
\end{align}
taking into account the boundedness of the symbol $B$ and hence of its pseudoinverse $B^\dagger$ (see \cite[Algorithm 1.3.1]{CM}).

Then, there exists a piecewise constant sequence $\lbrace h_n \rbrace$ with $h_n \to h$ in  $L^1(\Omega; \mathbb V)$  and $|Dh_n|(\Omega) \to |Dh|(\Omega).$  As  $\Omega$ is regular, we can extend $h$ outside $\Omega$ such that $|D h|(\partial \Omega)=0$ (see \cite[Theorem 2.16]{G}).  
Clearly, also $| \mathcal B h|(\partial \Omega) = 0,$ taking into account the characterization of the jump part, see \cite{ARS}.

Then, since $|Dh_n|(\Omega)$ is bounded, we have that $| \mathcal B h_n|(\Omega)= |B D h_n|(\Omega)$ is also bounded and therefore $ \mathcal B h_n \wsto  \mathcal Bh.$
We conclude that
\begin{align}\label{Ehn*}| \mathcal B h|(\Omega) \leq \liminf_{n \to +\infty} | \mathcal B h_n|(\Omega) \leq\limsup_{n\to +\infty} | \mathcal B h_n|(\Omega)\leq | \mathcal B h|(\Omega).
\end{align}
It is easily seen that $u_n \to g$ strongly in $L^1(\Omega;V)$ and, since $\mathcal B^a u_n =B\nabla u_n$, $\mathcal B^a h= B \nabla h$ and $\mathcal B^a h_n= B \nabla h_n=0$ it results  $\mathcal B^au_n= \mathcal B^a g+ \mathcal B^a h= G$. Taking into account the definition of $u_n$, we have
$\mathcal B u_n =B D u_n= \mathcal B g + \mathcal B h- \mathcal B h_n= B(D g+ D h-D h_n)$, thus, exploiting again the boundedness of the symbol $B$, the weak* convergence is immediately obtained. Concerning  the bound, in fact, we have
\begin{align*}
|\mathcal B u_n|(\Omega)\leq |\mathcal B g|(\Omega)+ |\mathcal B h|(\Omega)+ |\mathcal B h_n|(\Omega) \leq \|B\|( |D g|(\Omega)+ |Dh|(\Omega) +|D h_n|(\Omega)),
\end{align*}
where $\|B\|$ denotes the Frobenius norm of $B$.
 Up to the choice of a suitable tail of $(u_n)_n$, this  gives the estimate, in view of \eqref{Ehn*} and \eqref{817bis}.
\end{proof}

An analogous result holds when $\B$ is a second-order differential operator. In this case, the proof relies on a version of Theorem \ref{silhavy_app} in the space of functions with bounded Hessian $BH$. More precisely, an extension of Theorem \ref{Al} to the space $BH(\Omega)$ was established in \cite[Theorem 3.2]{FHP}, which we recall below. Recall also that the space $BH$ of functions with bounded Hessian is defined by
   \begin{equation}
     \notag
     BH(\Omega; V):=\{f \in W^{1,1}(\Omega;V): D^2 f \in \mathcal M(\Omega;V\otimes E_2(\mathbb R^N)\}.
 \end{equation}

\begin{Theorem}{\rm{(\cite[Theorem 3.2]{FHP})}}
\label{thm:FHP}
For every $(u,U)\in SD_{\B}(\Omega)$, there exists a sequence $(u_n)_n\subset SBH(\Omega;\mathbb{R}^d)$
such that $u_n\to u$ in $W^{1,1}(\Omega;\mathbb{R}^d)$ and $\nabla^2 u_n\stackrel{*}{\rightharpoonup} U$ in $\mathcal{M}(\Omega)$,
\[
\sup_n \|u_n\|_{BH}
\leq
C\left(\|u\|_{BH}+\|U\|_{L^1}\right),
\]
for some constant $C>0$.
\end{Theorem}

The counterpart of Theorem \ref{Al} for operator $\B$ of second order reads as follows. 

\begin{Theorem}\label{ourappthm:k=2}
 Let $\mathcal{B}$ be a $\mathbb C$-elliptic linear homogeneous differential operator of second order, whose symbol in jet form is $B^2$.  Let ${\rm Im}B^2\subseteq W$ be as in \eqref{ImBk}. Let  $(g,G)\in SBV^{\mathcal B}(\Omega)\times L^1(\Omega;{\rm Im} B^2)$.
 Then,  there exists a sequence $(g_n)_n \subset SBV^{\B}(\Omega)$ such that
$g_n \wsto g$ in $BV^{\mathcal B}(\Omega)$ and $\mathcal B^a g_n \wsto G$ in $\mathcal M(\Omega;{\rm Im}B^2)$.
\end{Theorem}
\begin{proof}[Proof]
Given $G \in L^1(\Omega;{\rm Im} B^2)$ and $g \in SBV^\B(\Omega)$, it results that $G(x)-\mathcal B^a g(x)= G(x)-B^2 \nabla^2 g(x) \in {\rm Im} B^2$ for a.e. $x \in \Omega$. Let $(B^2)^\dagger$ the Moore-Penrose pseudo inverse of $B^2$, then by Theorem \ref{thm:FHP} there exists $(f_n)_n \subset SBH(\Omega;V)$ such that $f_n \to 0$ in $W^{1,1}(\Omega ;V)$ and $\nabla^2 f_n \wsto (B^2)^\dagger(G(x)-\mathcal B^a g(x))$ in $\mathcal M(\Omega;V\otimes E_2(\mathbb R^N))$.
Consequently $\mathcal B^a f_n= B^2 \nabla^2 f_n\wsto G(x)-\mathcal B^a g(x)$ in $ \mathcal M(\Omega;{\rm Im}B^2)$. 
The sequence $g_n:= g+ f_n$ proves our claim.
\end{proof}

The Global Method and the approximation theorems above allow us to define the energy associated to a $\mathcal B$-structured deformation, i.e. to an element $(u, U) \in SD_\B(\Omega)$, as \eqref{def:SD}.
To this end, we start assuming that $\B$ is a $\mathbb C$-elliptic linear homogeneous, differential operator of order $k=1$ or $k=2$, as in \eqref{Bdef}.
We also need to  introduce the following notation: for $k=1,2$, we set $\hbox{grad}_{\B}^k u \coloneqq (\nabla^{k-1} u, \B^{a} u)$ and $j^ku\coloneqq ((\nabla^{k-2}u,\nabla^{k-1} u)^{+}, (\nabla^{k-1} u)^{-})$. Moreover, we introduce the vector space $\mathcal{Z}^k_{W}\coloneqq (\R^{d\times N})^{k-1}\times W$, i.e., for $k=1$ $\mathcal{Z}^1_{W}=W$ and for $k=2$ $\mathcal{Z}^2_{W}=\R^{d\times N}\times W$.

For every $\e>0$ we consider the energy functional $F_\e$  associated to a {\rm classical} deformation $u \in SBV^\B(A)$, given by
\begin{equation}
     \label{eq:functionalFei}
     F_\e(u; A) \coloneqq
     \begin{cases}
         \int_A f_\e(x, u, \hbox{grad}_{\B}^k u)dx + \int_{J_{\nabla^{k-1}u}\cap A} g_\e(x, j^ku, \nu_{\nabla^{k-1}u}) d{\mathcal H}^{N-1}, & u\in SBV^{\B}(\Omega);\\
            \\
            +\infty, & \hbox{otherwise},
     \end{cases}
 \end{equation}
where $f_\e$ and $g_\e$ satisfy the following assumptions
     \begin{itemize}
         \item[(T1)] $f_\e: \Omega \times \mathbb R^d \times \mathcal{Z}^k_{W}  \to [0, \infty)$ is measurable in $x$ and continuous in all other variables and there exists a constant $C>0$ such that
              \begin{equation*} 
                  \frac{1}{C} |\Lambda|\leq f_\e(x, a, \xi^{k-1}, \Lambda ) \leq C(1+|\Lambda|),
              \end{equation*}
        for all $(x, a, \xi^{k-1}, \Lambda ) \in\Omega \times \mathbb R^d\times \mathcal{Z}^k_{W}$;
        \item[(T2)] $g_\e: \Omega \times\mathbb R^d\times \mathbb R^{d\times N^{k-1}}\times\mathbb R^{d\times N^{k-1}}\times\mathbb S^{N-1} \to [0, \infty)$ is continuous and there exists a constant $C>0$ such that
           \begin{equation*} 
               \frac{1}{C} |\xi-\eta| \leq g_\varepsilon(x, a, \xi, \eta, \nu) \leq C(1+ |\xi-\eta|),
           \end{equation*}
         for all $(x, a, \xi, \eta,\nu) \in\Omega \times \mathbb R^d \times \mathbb R^{d\times N^{k-1}}\times \mathbb R^{d\times N^{k-1}}\times \mathbb S^{N-1}$. Moreover, $g_\e$ is symmetric, i.e., 
               \begin{equation}
                   \notag
                   g_\e(x, a, \xi, \eta, \nu) = g_\e(x, a, \eta, \xi, -\nu). 
               \end{equation}
         \end{itemize}
Here, we use the notation that $\nabla^0u=u$, and $\nabla^{-1} u=0$.

\begin{Remark}
    \begin{itemize}
        \item[{\rm (i)}] When $i=0$, the density $f_\e$ depends only on $(x, u, \B^a u)\in \Omega\times \R^d\times W$.
        \item[{\rm (ii)}] The family of functionals $F_\varepsilon^i$ encompasses two important examples:  the case in which the densities $f_\e$ and $g_\e$ are constant, i.e., $f_\e=f_0$ and $g_\e=g_0$, for some $f_0$ and $g_0$ satisfying Assumptions {\rm (T1)} and {\rm (T2)}, and the case where $f_\e$ and $g_\e$ are of the form 
            \begin{equation}\notag 
                f_\e(x,u, \text{grad}_\B^k u)\coloneqq f\left(x,\frac{x}{\e}, u, \text{grad}_\B^k u \right) \; \hbox{and}\; g_\e(x, j^ku)\coloneqq g\left(x, \frac{x}{\e}, j^ku\right), 
            \end{equation}
            with $f_\e$ and $g_\e$ being $[0,1)^N$-periodic in the second variable.
    \end{itemize}  
\end{Remark}

Now, we introduce the functional $\mathcal{F}_{SD_{\B}}:SD_{\B}(\Omega)\times\mathcal{A}(\Omega) \to[0,+\infty)$ given by
\begin{equation}\label{def:relFU}
\begin{aligned}
\mathcal{F}_{SD_{\B}}(u,U;A)
:=\inf\Big\{
&\liminf_{n\to+\infty} F_{\e}(u_\varepsilon;A):
u_\varepsilon\to u \text{ in }W^{k-1, 1}(\Omega;\mathbb{R}^d), \\
&\;\;\;\;\;\;\;\;\;\;\;\;\;\;\;\;\;\B u_\varepsilon\stackrel{*}{\rightharpoonup}U
\text{ in }\mathcal{M}(\Omega;W), k=1,\, 2\Big\},
\end{aligned}
\end{equation}
with $F_{\e}$ being as in \eqref{eq:functionalFei}.

In the case $f_\e= f_0$ and $g_\e=g_0$ the functional in \eqref{def:relFU} represents the energy associated to a $\mathcal B$-Structured Deformation, following the 
proposal in \cite{CF1997}. The idea consists in giving the structured deformation $(g, G)\in SD_\B$ the energy of the most economical
energetic approximation by means classical deformations, i.e. an asymptotic of the energies associated to the fields $u \in SBV^\B$, which approximate in the sense of either Theorem \ref{ourappthm} or Theorem \ref{ourappthm:k=2} the given $\B$ Structured Deformation $(g, G)$.

\begin{Theorem}
Let $\B$ be a homogeneous linear partial differential operator of order $k=1$ or $k=2$. Let $F_{\e}$ and $\mathcal{F}_{SD_\B}$ be the functionals defined by \eqref{eq:functionalFei} and \eqref{def:relFU}. Assume Assumptions {\rm (T1)} and {\rm (T2)}.  Then, there exist functions
\[
f:\Omega\times\mathbb{R}^d\times \mathcal{Z}_{\R^{d\times N}}^k
\to[0,\infty)
\]
and
\[
g:\Omega\times\mathbb{R}^d
\times\mathbb{R}^{d\times N^{k-1}}
\times\mathbb{R}^{d\times N^{k-1}}
\times S^{N-1}
\to[0,\infty)
\]
such that
\begin{equation}\label{eq:intrep-Def}
\begin{aligned}
\mathcal{F}_{SD_\B}(u,U;A)
&:=\int_A
f\bigl(x,u(x),\text{grad}^k_{\nabla^k},U\bigr)\,dx +
\int_{J_{\nabla^{K-1} u}\cap A}
g\bigl(x,j^ku, \nu_{\nabla^{k-1}u}\bigr)
\,d\mathcal{H}^{N-1}(x),
\end{aligned}
\end{equation}
for all $(u, U)\in SD_{\B}(\Omega)$
and $A\in\mathcal{A}(\Omega)$.
\end{Theorem}

\begin{proof}
    To show \eqref{eq:intrep-Def}, it suffices to prove that $\mathcal{F}_{SD_\B}$ satisties Assumptions ($\tilde H1$)-($\tilde H 4$) of Theorem \ref{thm:GMSDB}. Indeed,  the proof of Assumptions ($\tilde H2- \tilde H 4$) follows the same arguments of those deployed in Lemma \ref{lemma:lsc}, Lemma \ref{lemma:loc} and Lemma \ref{lemma:growth}, relying on assumptions (T1) and (T2). Note that in this setting we have to take into account not only the convergence of $u_\e \to u$ in $W^{k-1,1}$ but also $\mathcal B u_\e\wsto U$,
    (see also \cite[Lemmas 5.1, 5.2 and 5.3]{FHP} for an analogous result in $SD_2$).
   
     In order to prove that $\mathcal F_{SD_\B}$ satisfies Assumption ($\tilde H1$), we make use of the De Giorgi-Letta criterion (cf. Lemma \ref{DGL}). The only condition that one needs to prove is the nested subadditivity, i.e., for any open sets $A,B,C\in\mathcal{A}(\Omega)$ with $\overline{A}\subset B\subset C$ we have
         \begin{equation}
             \notag
              \mathcal{F}_{SD_\B}(u, U;C)\leq \mathcal{F}_{SD_\B}(u, U;C\setminus\overline{A})+\mathcal{F}_{SD_\B}(u, U;B).
         \end{equation}
         To this end, it is sufficient to argue as in Lemma \ref{lemma:growth}, taking into account the presence of the field $U$.
\end{proof}

\begin{Examples}
\label{examples}
   Several mechanical models are covered by our analysis thus providing a unified theory of Structured Deformations $SD_\B$. In particular the case when $\mathcal B$ is a homogeneous, $\mathbb C$-elliptic, linear differential operator of order $1$ covers \begin{itemize}
   \item[{\rm (a)}] for $\mathcal B = D$, the classical theory of structured deformations in {\rm \cite{DPO, CF1997}};
   \item[{\rm (b)}] for $\mathcal B= \mathscr{E}$ (where the latter represents the distributional symmetrized gradient) the theory of Structured Deformations in linear elasticity presented in {\rm \cite{FMZ}};
   \item[{\rm (c)}] when $\mathcal B ={\rm dev} \mathscr{E}= \mathscr{E}- \tfrac{{\rm div}}{N}Id$ a theory of structured deformations for the deviatoric operator, see {\rm \cite{CVG}} for energetic models related to soap films and rubber like materials.
   \end{itemize}
   The case of the $\mathbb C$-elliptic operator $\mathcal B$ of order two, allows to deal with the theories of nonsimple materials {\rm \cite{T1,T2}}, where also gradient disarrangements play a role, as suggested in {\rm \cite{O17}}. In particular, in \eqref{eq:functionalFei}, one can consider bulk free energies densities $f_\e(x, u(x), {\rm grad}^2_{\B} u(x))$ of the form $f_0(x, \nabla u(x), \B^{a}u(x))$ as suggested in {\rm \cite[eq. (97)]{O17}}. To fix ideas, one could consider several models as recalled in {\rm \cite[Sections 2, 5]{B}} which allow the treatment of various types of second gradient materials, e.g. materials, as in {\rm \cite[Section 3.10]{O17}}, whose deformations are ($\B$)-"submacroscopically affine". For  $\mathcal B= D^2$ one recovers the, by now, classical theory of second order structured deformations as in {\rm \cite{OP, FHP}}. 
   
   Other choices for the dependence on the Hessian tensor $D^2$ are also possible in the case of small deformations. In particular, one can consider the models proposed by {\rm \cite{M64,M65,ME68}} as recalled in {\rm \cite{B}}:
   \begin{itemize}
   \item[{\rm (i)}] $\B= D^2$, {\rm second displacement (distributional) gradient} which encodes the dependence on the entire gradient, as in the case of large deformations, already mentioned above;
   \item[{\rm (ii)}] $\mathcal B = D (\mathscr{E})$, {\rm (distributional) gradient of the (distributional) strain tensor}, which in the end is another example of $\mathbb C$-elliptic operator written in jet form as $\mathcal B= B^2_1 D^2$ with ${B^2_1}$ of maximal rank. In particular this second case presents only a rearrangement of the hessian tensor $D^2$;
   \item[{\rm (iii)}]  $\B = \mathscr{E}(\mathscr{E})$, (computed row-wise) which is the {\rm completely symmetric part of the (distributional) gradient of the (distributional) strain tensor}
in this case $\mathcal B= B^2_2 D^2$ with $B^2_2$ of lower rank, but still providing a $\mathbb C$-elliptic operator. This latter operator also has the name of completely symmetrized Hessian.
   \end{itemize}
Note that in {\rm (i)} and {\rm (ii)}, under standard regularity assumptions on $\Omega$ and due to \eqref{BVBequiv}, the spaces $\{u\in LD(\Omega)\; : \; D^2u\in \mathcal{M}(\Omega; \mathbb R^N\otimes E_2(\R^N))\}$ and $\{u\in LD(\Omega)\; : \; D(\mathscr{E})u\in \mathcal{M}(\Omega; \mathbb R^N\otimes E_2(\R^N))\}$ coincide with the space $BH(\Omega)$ (recall that $LD(\Omega):=\{u \in L^1(\Omega;\mathbb R^N), \mathscr{E}u \in L^1(\Omega;\mathbb R^{N \times N})\}$, see {\rm \cite{T} for details}).

While ${\rm (iii)}$, arguing as above, $\{u \in LD(\Omega): \mathscr{E}(\mathscr{E} u) \in \mathcal M(\Omega; E_3(\mathbb R^N)\}= BV^{\mathscr{E}(\mathscr{E})}(\Omega)$, cf. { \rm \cite[Example 8.4]{ARS} and \cite{VS}}. 
   
We conclude observing that
    our results are applicable to families of functionals. This latter approaches proves to be useful when dealing with the study of the asymptotic behaviour of materials with finely distributed heterogeneities, i.e. in particular it allows to extend the homogenization results in {\rm \cite{AMMZ}, \cite{BMZ}} to more general Structured Deformations.
\end{Examples}

\section{Appendix} \label{sec:appendix}

The following is the variant of the De Giorgi Letta criterion as quoted in \cite[Theorem 1.53]{AFP} exploited in Section \ref{sec:applications}.
\begin{lemma}\label{DGL} Let $\mathcal{A}(\Omega)$ be the family of open subsets of $\Omega$ and let \(\lambda : \mathcal{A}(\Omega) \to [0,\infty) \) be an increasing set function such that: \begin{itemize} \item[(i)] for all $A,B,C \in \mathcal{A}(\Omega)$ with \( A \subset B \subset C, \) there holds \[ \lambda(C) \leq \lambda(C \setminus A) + \lambda(B); \] \item[(ii)] \[ \lambda(A \cup B)=\lambda(A)+\lambda(B), \] for all $A,B \in \mathcal{A}(\Omega)$ such that \(A \cap B=\varnothing; \) \item[(iii)] there exists a measure \[ \mu : \mathcal{B}(\Omega)\to [0,\infty) \] such that \(\lambda(A)\leq \mu(A) \) for all $A\in\mathcal{A}(\Omega)$, where $\mathcal{B}(\Omega)$ denotes the family of Borel subsets of $\Omega$. \end{itemize} Then $\lambda$ is the restriction to $\mathcal{A}(\Omega)$ of a measure defined on $\mathcal{B}(\Omega)$. \end{lemma}

In particular, with an abuse of notation $\lambda$ can be extended to $\mathcal B(\Omega)$, as follows
\[
\lambda(B) := \inf\{\lambda(A) : A \in \mathcal{A},\ A \supset B\},
\]
and it is an outer measure. The restriction of $\lambda$ to the Borel sets of $X$, $\mathcal{B}(X)$, 
is a positive measure.

The following result, whose proof can be found in \cite{FMARMA} is recalled for the readers' convenience and it will be used in the sequel
\begin{lemma}\label{FM}
Let $\lambda$ be a non-negative Radon measure in $\mathbb R^N$. For $\lambda$ a.e. $x_0 \in \mathbb R^N$, for every $0 < \delta < 1$ and for every $\nu \in \mathbb S^{N-1}$, the following holds
$$\limsup_{\eps \to 0^+}\frac{\lambda(Q_\nu(x_0,\delta \eps))}{\lambda(Q_\nu(x_0,\eps))} \geq \delta^N,$$
so that
$$\lim_{\delta \to 1^-}\limsup_{\eps \to 0^+}\frac{\lambda(Q_\nu(x_0,\delta \eps))}{\lambda(Q_\nu(x_0,\eps))} = 1.$$
\end{lemma}

Formulas \eqref{jumppt}-\eqref{jumpptk-2} can be deduced by \cite[Theorem 2.4, Remark 2.3, Corollary 2.8, and Lemma 2.17]{ARS}. However, we present a more explicit proof, whose strategy is based on \cite[Exercise 3.16]{AFP}.

\begin{Proposition}\label{prop:estAFP}
 Let $\B$ be a differential operator, $\mathbb{C}$-elliptic and of order $k$, with $k>1$. Let $u\in BV^{\B}(\Omega)$. Then, for $\mathcal{H}^{N-1}$-a.e. $x\in J_{\nabla^{k-1}u}$, we have 
         \begin{equation}
             \label{eq:propjumpt}
             \lim_{r\to 0} \frac{1}{\rho^N}\int_{B^{\pm}_{\rho}(x, \nu)} \left(\frac{|u(y)-P^{k-2}_x(y-x)-\frac{1}{(k-1)!}\langle{(\nabla^{k-1})}^{\pm}(x), (y-x)^{\otimes^{k-1}}\rangle|}{|y-x|^{k-1}}\right)^{1^{*}}dy=0,
         \end{equation}
    where $P^{k-2}_x(y-x)$ is the polynomial defined as in \eqref{Pkx0def} and the half balls $B^{\pm}_{\rho}(x, \nu)$ are defined in \eqref{def:halfball}.
   .
\end{Proposition}
\begin{proof} 
Since the proof of \eqref{eq:propjumpt} is iterative, we first prove \eqref{eq:propjumpt} is for an operator $\B$ of the second order and then for any order $k>2$.
\smallskip

{\bf Step 1.} In this step, we are going to prove the following
        \begin{equation}
             \notag
             \lim_{r\to 0} \frac{1}{\rho^N}\int_{B^{\pm}_{\rho}(x, \nu)} \left(\frac{|u(y)-u(x)-\langle{\nabla}^{\pm}u(x), y-x\rangle|}{|y-x|}\right)^{1^{*}}dy=0,
         \end{equation}
which is the analogous of \eqref{eq:propjumpt}  for operators $\B$ of order $k=2$.  We prove the equality just for $B^{+}_{\rho}(x, \nu)$, being the other equality completely analogous.   First, note that, for $\mathcal{H}^{N-1}$-a.e. $x\in \Omega$, $B^{\pm}_{\rho}(x, \nu)= B_{\rho}(x, \nu)\cap H^{+}$, where $H^{+}\coloneqq\{y\in \Omega : (y-x)\cdot \nu>0 \}$ and $B_\rho(x)$ can be written as follows 
\[
 B_\rho(x)=\bigcup_{i=0}^{\infty}
 \left(B_{\rho/2^i}\setminus B_{\rho/2^{i+1}}\right)(x),
\]
which implies that is $ y\in\left(B_{\rho/2^i}\setminus B_{\rho/2^{i+1}}\right)(x)$ then $|y-x|>\frac{\rho}{2^{i+1}}.$ Therefore,
\begin{align*}
&\int_{B^{\pm}_{\rho}(x, \nu)}
 \left(\frac{\bigl|u(y)-u(x)-\langle\nabla^+u(x),y-x\rangle\bigr|}
 {|y-x|}\right)^{1^{*}}\,dy\\
&=\sum_{i=0}^{\infty}
 \int_{(B_{\rho/2^i}\setminus B_{\rho/2^{i+1}})(x)\cap H^{+}}
 \left(\frac{\bigl|u(y)-u(x)-\langle\nabla^+u(x),y-x\rangle\bigr|}
 {|y-x|}\right)^{1^{*}}\,dy\\
&\leq \sum_{i=0}^{\infty}
 \frac{1}{(\rho/2^{i+1})^{1^{*}}}
 \int_{(B_{\rho/2^i}\setminus B_{\rho/2^{i+1}})(x)\cap H^{+}}
 \bigl|u(y)-u(x)-\langle\nabla^+u(x),y-x\rangle\bigr|^{1^{*}}\,dy.
\end{align*}
Hence,
\begin{align*}
&\frac1{\rho^N}\int_{B^{\pm}_{\rho}(x, \nu)}
 \left(\frac{\bigl|u(y)-u(x)-\langle\nabla^+u(x),y-x\rangle\bigr|}
 {|y-x|}\right)^{1^{*}}\,dy\\
&\leq 2^{1^{*}}\sum_{i=0}^{\infty}2^{-iN}
 \frac{1}{(\rho/2^i)^{1^{*}+N}}
 \int_{(B_{\rho/2^i}\setminus B_{\rho/2^{i+1}})(x)\cap H^{+}}
 \bigl|u(y)-u(x)-\langle\nabla^+u(x),y-x\rangle\bigr|^{1^{*}}\,dy\\
&\leq 2^{1^{*}}\sum_{i=0}^{\infty}2^{-iN}
 \frac{1}{(\rho/2^i)^{1^{*}+N}}
 \int_{B_{\rho/2^i}(x)\cap H^{+}}
 \bigl|u(y)-u(x)-\langle\nabla^+u(x),y-x\rangle\bigr|^{1^{*}}\,dy.
\end{align*}
Note that due to the embedding theorems \cite{GR19},  $u\in BV^{\B}(\Omega)$ implies that  $u\in W^{2,1^{\ast}}(\Omega; W)$. In particular, $\nabla u\in L^1(\Omega; W)$. Therefore, making use of the Poincar\'{e} inequality in the Sobolev spaces along with the results in \cite{D} (see also \cite[Theorem 2.4]{ARS}), we deduce that  
\begin{align*}
&
 \int_{B_{\rho/2^i}(x)\cap H^{+}}
 \bigl|u(y)-u(x)-\langle\nabla^+u(x),y-x\rangle\bigr|^{1^{*}}\,dy
 \\
&\qquad\leq C\rho^{1^\ast}\,|\nabla u-\nabla^+u(x)\mathcal L^N|
 \bigl(B_{\rho/2^i}(x)\cap H^{+}\bigr)
=\rho^{1^\ast}o\!\left((\rho/2^i)^N\right).
\end{align*}
Then, for $\mathcal H^{N-1}$-a.e. $x\in\Omega$,
\[
 \lim_{\rho\to0}\frac1{\rho^N}
 \int_{B^{\pm}_{\rho}(x, \nu)}
 \frac{\bigl|u(y)-u(x)-\langle\nabla^+u(x),y-x\rangle\bigr|^{1^{*}}}
 {|y-x|^{1^{*}}}\,dy=0,
\]
concluding the proof.
\medskip

{\bf Step 2.} For operators $\B$ of order $k>2$, \eqref{eq:propjumpt} follows by applying iteratively the Poincar\'{e} ineqality along with the fact that $u\in W^{k-1, 1^{\ast}}(\Omega; W)$.    
\end{proof}

An application of the H\"older inequality yields
\begin{Corollary}\label{cor:jumppt}
    Under the same assumptions of Proposition \ref{prop:estAFP}, we also have that 
             \begin{equation}
             \notag
             \lim_{r\to 0} \frac{1}{\rho^N}\frac{1}{\rho^{k-1}}\int_{B^{\pm}_{\rho}(x, \nu)} |u(y)-P^{k-2}_x(y)-\frac{1}{(k-1)!}\langle{(\nabla^{k-1})}^{\pm}(x), (y-x)^{\otimes^{k-1}}\rangle|dy=0,
         \end{equation}
with $B^{\pm}_{\rho}(x, \nu)$ being defined as in \eqref{def:halfball}.
\end{Corollary}

\color{black}

\subsection*{Acknowledgements}
The authors thank D. R. Owen for his deep insights into the theory of structured deformations and A. Arroyo-Rabasa, M. Caroccia and N. Van Goethem for fruitful discussions on the topics addressed in this manuscript.
L. D'Elia acknowledges support of the Austrian Science Fund (FWF) project 10.55776/ESP1887024. She also acknowledge the support of Gruppo Nazionale per l'Analisi Matematica, la Probabilit\'a e le loro Applicazioni-INdAM through Project Professori Visitatori, CUP E53C25002010001.
E. Zappale is a member of GNAMPA-INdAM, whose support is gratefully acknowledged. She is very grateful to TU Wien for its hospitality and for the support through project 10.55776/ESP1887024.

\end{document}